\documentclass[journal,onecolumn]{IEEEtran}

\usepackage{float}
\usepackage{units}
\usepackage{amsmath}
\usepackage{graphicx}

\usepackage[unicode=true,
 bookmarks=true,bookmarksnumbered=true,bookmarksopen=true,bookmarksopenlevel=1,
 breaklinks=false,pdfborder={0 0 0},pdfborderstyle={},backref=false,colorlinks=false]
 {hyperref}

\usepackage{placeins}

\begin{document}

\title{A Novel Consensus-based Distributed Algorithm for Economic Dispatch
Based on Local Estimation of Power Mismatch }

\author{Hajir~Pourbabak,~\IEEEmembership{Student Member,~IEEE,} Jingwei
Luo,~\IEEEmembership{Student Member,~IEEE,} Tao~Chen,~\IEEEmembership{Student Member,~IEEE}
and~Wencong~Su,~\IEEEmembership{Member,~IEEE}\thanks{Hajir~Pourbabak, Jingwei Luo, Tao Chen and Wencong Su are with the
Department of Electrical and Computer Engineering, College of Engineering
\& Computer Science, University of Michigan-Dearborn, MI, USA, e-mail:
\protect\href{http://hpourbab@umich.edu}{hpourbab@umich.edu}}}
\maketitle
\begin{abstract}
\textbf{This paper proposes a novel consensus-based distributed
control algorithm for solving the economic dispatch problem of distributed
generators. A legacy central controller can be eliminated in order
to avoid a single point of failure, relieve computational burden,
maintain data privacy, and support plug-and-play functionalities.
The optimal economic dispatch is achieved by allowing the iterative
coordination of local agents (consumers and distributed generators).
As coordination information, the local estimation of power mismatch
is shared among distributed generators through communication networks
and does not contain any private information, ultimately contributing
to a fair electricity market. Additionally, the proposed distributed
algorithm is particularly designed for easy implementation and configuration
of a large number of agents in which the distributed decision making
can be implemented in a simple proportional-integral (PI) or integral
(I) controller. In MATLAB/Simulink simulation, the accuracy of the
proposed distributed algorithm is demonstrated in a 29-node system
in comparison with the centralized algorithm. Scalability and a fast
convergence rate are also demonstrated in a 1400-node case study.
Further, the experimental test demonstrates the practical performance
of the proposed distributed algorithm using the VOLTTRON\texttrademark{}
platform and a cluster of low-cost credit-card-size single-board PCs.}
\end{abstract}




\section{Introduction}

Future power systems are equipped with a great
number of distributed generators (DGs), distributed energy storage
devices, dispatchable loads and advanced communication networks, which
increases the customer participation in the electricity market. As
a result, the optimal economic dispatch (ED) of future power systems
is becoming much more challenging. A sophisticated control is needed
to fully address the increasing customer participation at the edge
of the electric power system and the inability of existing practices
to accommodate these changes \cite{Fang2012,Potter2009,NationalAcademiesofSciencesandMedicine2016,Asadinejad}.

\begin{figure*}[tbh]
\includegraphics[width=0.8\paperwidth]{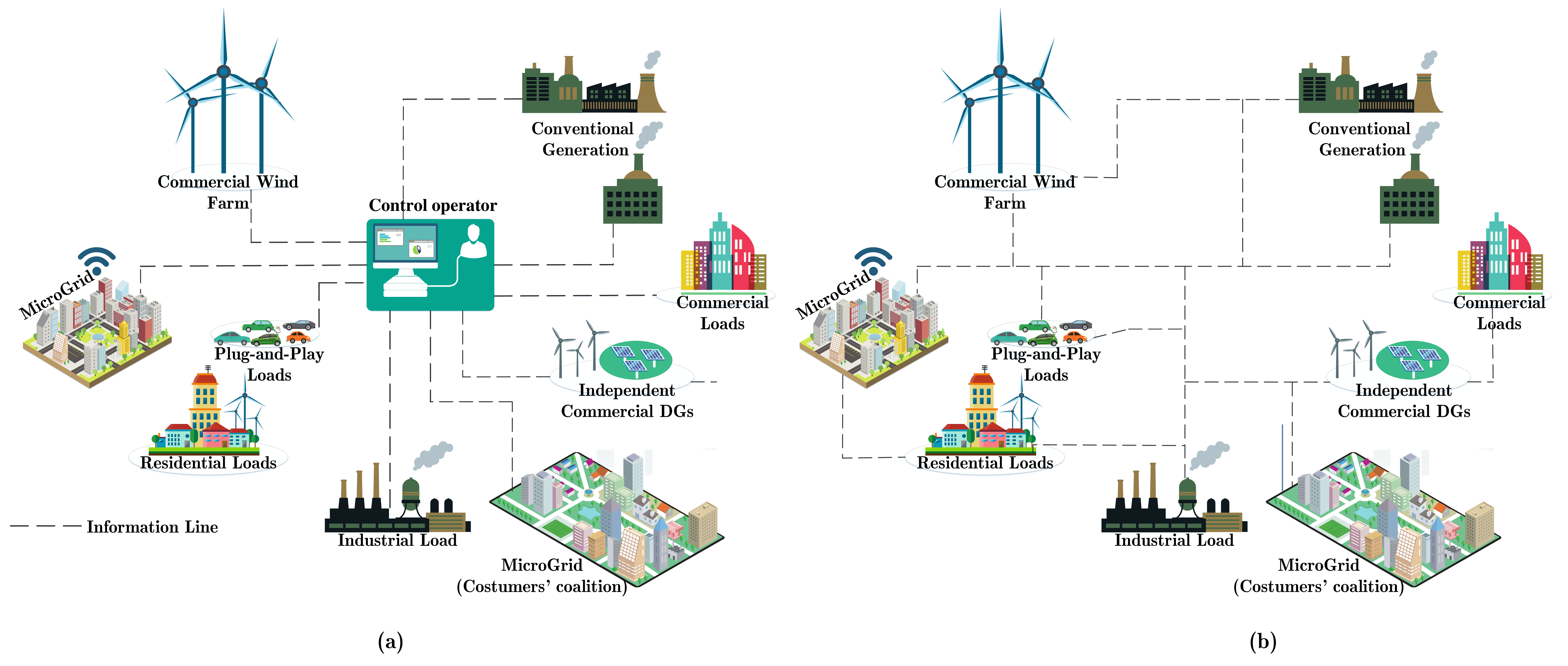}

\caption{Different approaches for power grids. (a) centralized approach, (b)
distributed approach (peer-to-peer connections)\label{fig:FuturePowerGrid}}
\end{figure*}

In a centralized ED, all participants (DGs and consumers) must release
their information to the central controller. As the market penetration
of DGs and consumers is continuously growing, the centralized algorithms
are no longer suitable due to the heavy computational burden \cite{Mudumbai2012c}.
Moreover, the centralized algorithms are not designed to support plug-and-play
functionalities of a large number of participants.

Consensus-based distributed approaches have been found to be practical
in many multi-agent applications, such as industrial systems, automated
high way systems, and computer networks \cite{Saber2003}. Consensus-based
distributed approaches are to find the global optimal decision by
allowing local agents to iteratively share information through two-way
communication links. All agents reach a consensus when they agree
upon the value of the information state \cite{Ren,Olfati-Saber2007,Binetti2014,Xu2017,Xu2016Su}.
The information state can be physical quantities or control signals
such as voltage, frequency, output power, incremental cost, and estimated
power mismatches \cite{Zhang2012d}. Figure \ref{fig:FuturePowerGrid}
compares the centralized and distributed methods for solving the ED
problem in power systems \cite{Pourbabakbook}. The major advantages
of distributed methods are summarized as follows:

\textbf{\textit{Scalability and Interoperability}}: As the number
of agents increases to hundreds of thousands, the legacy centralized
method faces certain challenges such as computational burden. As more
DGs are integrated into power systems, the centralized methods are
not suitable for such heterogeneous systems \cite{Zhang2012d,Zeng2014}.

\textbf{\textit{Monopoly and Monopsony}}: A central organization usually
has a kind of monopoly over the consumers and a monopsony toward the
electrical energy DGs. However, as it is discussed in various research
papers and academic notes, consumers have a much more apathetic role
than DGs do in electrical energy markets under centralized control
\cite{Kirschen2004,Samadi2012}.

\textbf{\textit{Privacy and Stealth Protection}}: One of the essential
and notable features of every competitive system is the equal opportunity
for competition for all players. Therefore, it should be ensured that
no private information is released by a third-party and no preferences
are given to one or more DGs in a competitive market \cite{Fang2012}.

\textbf{\textit{Computational Cost: }}Heavy computational load is
imposed on the central controller when dealing with a large multi-agent
network \cite{Zhang2012d,Xu2015b}. However, the computational load
can be shared among agents using distributed approaches.

\textbf{\textit{Single Point of Failure}}: A distributed approach
is robust to the single point of failure because there is no need
for a center for the supervisory of the entire system \cite{Zeng2014,Soediono1989}.

\textbf{\textit{Network Topology: }}In future power systems, the communication
network and power network are subject to frequent change \cite{Kazemi2016}.
Accordingly, there are serious doubts about the ability of centralized
methods to handle the variable topology.

The literature review shows a growing interest in distributed algorithms
in the field of power systems. A distributed algorithm
for solving the ED problem that considered thermal generation and
random wind power was discussed in \cite{Guo2014} for a smart grid.
A distributed optimization problem with local constraints was studied
in \cite{Yi2015} to obtain the optimal solution for a load sharing
problem. Hug and Kar \cite{Hug2015a} introduced a consensus-based
distributed energy management approach to handle loads of a micro-grid.
Rahbari-Asr \textit{{et al}}.
investigated an incremental welfare consensus-based algorithm that
manages both loads and DGs \cite{Rahbari-Asr2014a}. In addition,
they developed a cooperative distributed demand management system
in \cite{Rahbari-Asr2014} based on Karush\textendash Kuhn\textendash Tucker
conditions for plug-in electric vehicle charging. A distributed optimal
power flow was introduced for a large-size distribution system in
\cite{DallAnese2013}. Mudumbai \textit{{et al}}.
\cite{Mudumbai2012c} worked on the distributed algorithm for both
optimal ED and frequency control with a large number of intermittent
energy sources. The authors in \cite{Multibuyer2015} proposed a distributed
real-time demand response for a multi-agents system to maximize social
welfare. They defined sub-problems by decomposing the main optimization
problem and locally solving each sub-problem. Binetti \textit{{et
al}}. \cite{Binetti2014} presented a distributed
method to solve the ED problem considering power losses of transmission
lines. Two consensus protocols were run in parallel to reach the consensus
and estimate the power-mismatch considering the lines' losses. Elsayed
and El-Saadany proposed a fully decentralized method that can solve
both the convex and the practical non-convex ED problem. They also
considered transmission losses in \cite{Elsayed2015}. In \cite{Zhang2015b},
a distributed constrained gradient approach was proposed to consider
both the equality and inequality constraints for online optimal generation
control. Another algorithm has been introduced in \cite{Yang2013a}
to help energy producers to collect information regarding power mismatch
between generation and consumption. This mechanism finally converges
incremental cost to an optimal value. N. Cai \textit{{et al.}} \cite{Cai2012} offers a decentralized approach
for the economic dispatch problem of a micro-grid in which, various
agents use local information or receive some information from their
immediate neighbors. However, this paper does not cover the power
balance within micro-grid because it is assumed that power balance
has already been met. One of the distributed techniques for solving
non-convex economic dispatch problem has been introduced in \cite{Binetti2014b}.
In this technique, a leader-less consensus based distributed algorithm
share bids among different agents based on an auction mechanism. G.
Binetti \textit{et al}. considered
transmission losses, valve-point loading effect, multiple fuel, prohibited
operating zones in non-convex ED problem. A primal-dual perturbation
method is proposed by T.H Chang \textit{{et al}}.
\cite{Chang2014} enabling multi-agents system to reach a global optimum.
In this method, all agents try to estimate functions of global cost
and constraints based on distributed approach.

\begin{figure}[tbh] \centering
\includegraphics[width=0.5\columnwidth, ]{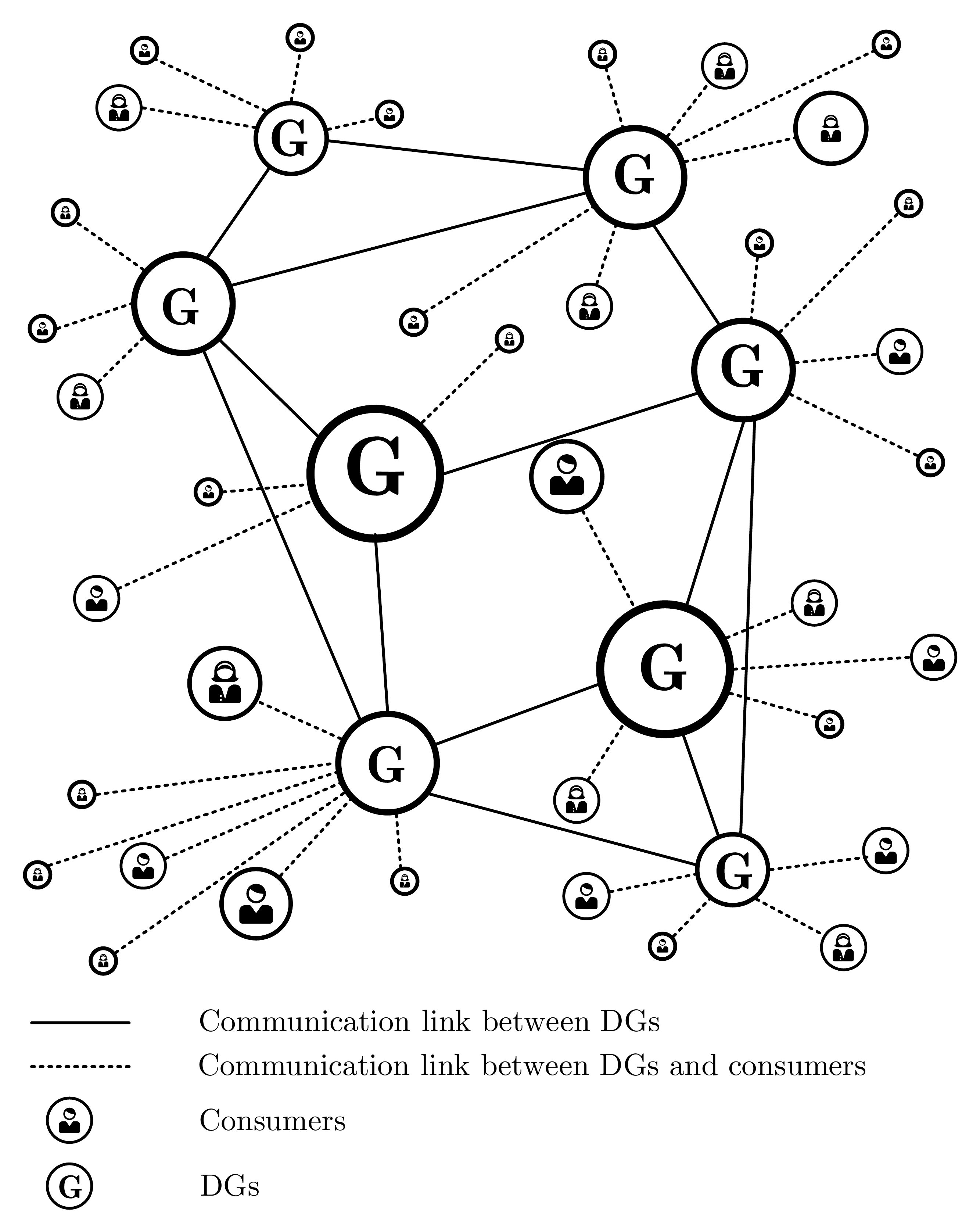}

\caption{The illustrated communication network for the proposed distributed
method.\label{fig:Theconceptofproposedmethod}}
\end{figure}

From the control algorithm\textquoteright s point of view, one common
issue with existing distributed consensus algorithms is their slow
convergence (high iteration numbers), privacy, high connectivity requirement
and complexity. For instance, in most of previous works the incremental
cost ($\lambda$) and output power are shared \cite{Multibuyer2015},
which violate privacy principles. While agents in a multi-agent system
may not share private information with center (third-party), they
share the private information with other agents and it could be even
worse than the centralized method. In addition, distributed methods
cannot attract serious attentions in case of implementation if they
suffer from high connectivity, slow convergence and complexity. To
address the aforementioned limitations of existing distributed methods,
we propose a novel consensus-based distributed algorithm to maintain
data privacy and reduce computational time while solving for a large-scale
ED problem. The proposed approach need to share minimum information
(Power mismatch) among different agents in a multi-agent network,
thus, the privacy of all agents will definitely be improved. In other
words, the parameter of cost function, utility function, incremental
cost ($\lambda$) and output power etc. will not be shared among
agents and, also, third-party are not at all able to access to these
parameters. In addition, computational cost will be decreased that
make the scalability possible for enormous multi-agent system. Figure
\ref{fig:Theconceptofproposedmethod} shows a general view of the
communication network for the proposed method. The main features of
the proposed distributed algorithm are as follows:
\begin{enumerate}
\item \textbf{\textit{Accuracy}}: The ED problem is solved in a fully distributed
manner without relying on a central controller. The solution accuracy
is validated by the benchmark results achieved by traditional centralized
methods.
\item \textbf{\textit{Privacy}}: DGs and consumers do not need to disclose
any private information (e.g., cost functions and utility functions)
with others. The estimated power mismatch between the total generation
and the total demand is the only information to be shared among all
DGs. In addition, consumers are not required to have a communication
channel between one another or to exchange information with all DGs.
They only need to connect to local DGs.
\item \textbf{\textit{Fast Convergence}}: The proposed algorithm outperforms
some existing distributed methods in terms of number of iteration
and computational time.
\item \textbf{\textit{Scalability}}:
The proposed distributed algorithm is particularly suitable for solving
large-scale optimization problems (e.g., >1,000 agents) within a short
period of time.
\item \textbf{\textit{Easy Implementation}}: This salient feature makes
it possible to deploy the proposed distributed methods in the field
at scale.
\end{enumerate}

The structure for the rest of this paper is organized as follows:
Section \ref{sec:SystemModelingandproblem} formulates the ED problem
as a global objective function, considering cost functions, utility
functions and constraints. Section \ref{sec:DistributedAlgorithmforED}
discusses graph theory, consensus-based distributed protocols, optimality
and convergence analysis of proposed algorithm . Section \ref{sec:PerformanceAssessment}
evaluates the solution performance using software simulation and experimental
testing. Section \ref{sec:ConclusionsandFuture} summarizes this paper
and presents the concluding remarks.

\section{System Modeling and Problem Description\label{sec:SystemModelingandproblem}}

The ED problem is a short-term resource allocation of a number of
DGs to meet the load requirement in a most cost-effective way. The
utility function of consumers, the cost function of DGs and their
surplus function are defined in this section. Then, the overall optimization
problem is formulated based on the defined model of economic players
in an electricity market.

\subsection{Utility Function, Marginal Benefit and Consumer's Surplus}

Each consumer in an electricity market has its own preferences for
energy consumption during different times of a day. These preferences
cause various levels of requested demand within an operation period.
There are some factors affecting the preferences of consumers in an
electricity market such as the instantaneous or average price of electricity,
temperature changes, the type of user and comfort level. The different
demand levels requested by consumers in response to these diverse
factors can be modeled by utility functions for mathematical purposes
\cite{Samadi2012,Multibuyer2015,Pourbabakconf}. In other words, the
utility function measures the satisfaction level or welfare of consumers
as a function of different types of performance (\textit{i.e.,} demand
level) to represent a consumer\textquoteright s preferences. In a
typical electricity market, the utility function $U_{j}(P_{j}^{Load})$
shows the level of satisfaction of the $j-th$ energy consumer, where
$P_{j}^{Load}$ is the demand of $j-th$ consumer.

Marginal benefit is an additional utility that an electrical consumer
will gain by getting one more unit ($\unit{MW}$) of electrical energy.
The consumption of more power increases the utility function if the
marginal benefit has a positive value; however; the consumption of
more power decreases the level of satisfaction if the marginal benefit
has a negative value \cite{Fahrioglu2001}. The common utility functions
are non-decreasing functions; thus, the marginal benefit is a non-negative
function. It means \cite{Samadi2010,Mohajeryami2016}:
\begin{equation}
\frac{\partial U_{j}\left(P_{j}^{Load}\right)}{\partial P_{j}^{Load}}\geq0,\qquad j\epsilon S_{D}\label{eq:nondecreasing}
\end{equation}

It is evidently proven that the marginal benefit of electrical consumers
is a non-increasing function because the marginal benefit will normally
decrease as users consume more energy. Thus, we have \cite{Samadi2010,Mohajeryami2016}:

\begin{equation}
\frac{\partial^{2}U_{j}\left(P_{j}^{Load}\right)}{\partial\left(P_{j}^{Load}\right)^{2}}\leq0,\qquad j\epsilon S_{D}\label{eq:nonincreasing}
\end{equation}

There are various types of utility function for single and multiple
goods, such as: \textit{Cobb-Douglas Utility Function}, \textit{Perfect
Substitutes}, \textit{Perfect Complements} and \textit{Quasilinear}
\cite{Nechyba2016}. In this paper, a quadratic mathematical model
satisfying (\ref{eq:nondecreasing}) and (\ref{eq:nonincreasing})
is used in Equation (\ref{eq:utilityfunction}) as a utility function
of the consumers. This utility function is customized for different
consumers based on parameters $\unitfrac[b]{\$}{kWh^{2}}$ and $\unitfrac[\omega]{\$}{kWh}$.
The larger $\omega$ is, the higher the utility function \cite{Multibuyer2015}.

\begin{gather}
U_{j}\left(P_{j}^{Load}\right)=\begin{cases}
\omega_{j}P_{j}^{Load}-b_{j}(P_{j}^{Load})^{2} & P_{j}^{Load}\leq\nicefrac{\omega_{j}}{2b_{j}}\\
\frac{\omega_{j}^{2}}{4b_{j}} & P_{j}^{Load}\geq\nicefrac{\omega_{j}}{2b_{j}}
\end{cases}\label{eq:utilityfunction}
\end{gather}

Consumer's surplus measures the welfare on the consumers' side; hence,
it is a measurement of the benefit, derived from the electricity market,
of an economic player on the consumption side \cite{Rahbari-Asr2014a}.
A consumer's surplus will be represented by (\ref{eq:consumersurplus}),
if the $j-th$ consumer pays $\lambda\;\unitfrac{\$}{kWh}$ for $P_{j}^{Load}\;\unit{kW}$
of the electrical energy. 

\begin{equation}
\mathcal{S}_{j}^{Load}\left(P_{j}^{Load}\right)=U_{j}\left(P_{j}^{Load}\right)-\lambda P_{j}^{Load},\qquad j\epsilon S_{D}\label{eq:consumersurplus}
\end{equation}

\begin{figure}[tbh] \centering
\includegraphics[width=0.7\columnwidth]{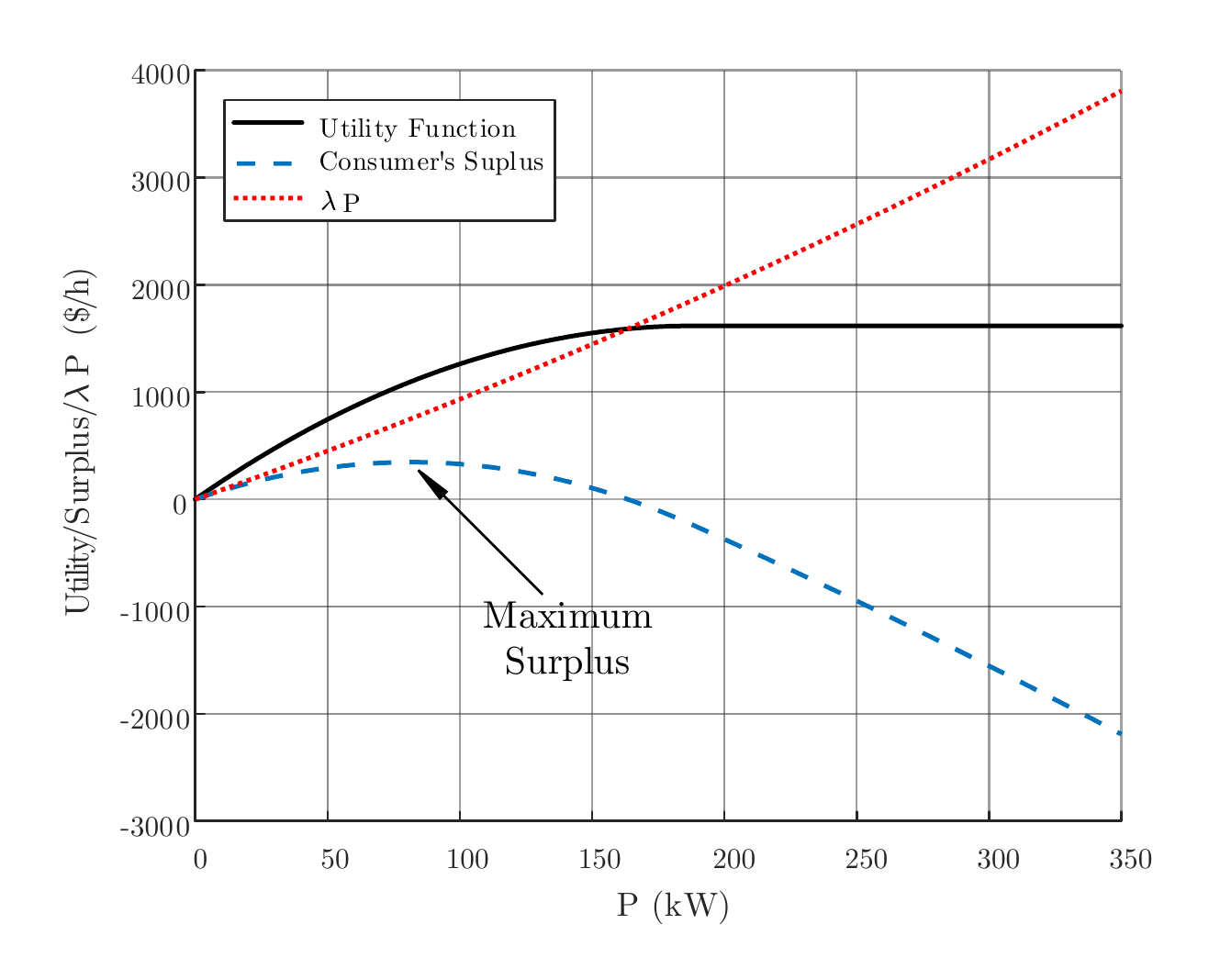}

\caption{Consumer's surplus and utility function\label{fig:Consumrsuplusutility}}
\end{figure}

Consumers attempt to maximize their own welfare in the market. Therefore,
they consume power at the maximum value of their concave surplus function.
Figure \ref{fig:Consumrsuplusutility} shows the graphical representation
of equations (\ref{eq:utilityfunction}) and (\ref{eq:consumersurplus}).

\subsection{Cost Function and DG's Surplus}

Generally speaking, a multiple piecewise linear or quadratic function,
known as cost function, is used to estimate the total cost of the
output-power of the energy providers, such as DGs. Using a proper
cost function is the best way to pre-evaluate the performance of a
DG and to solve an ED problem \cite{Kirschen2004,Pourbabakconf}.
Here, we consider a quadratic mathematical equation to model a typical
DG. $\unitfrac[\alpha]{\$}{kWh^{2}}$, $\unitfrac[\beta]{\$}{kWh}$
and $\unitfrac[\gamma]{\$}{h}$ are coefficients that customize the
cost function for each DG and $P_{i}^{Gen}$ is power generated by
the $i-th$ DG. 

\begin{equation}
C_{i}\left(P_{i}^{Gen}\right)=\alpha_{i}(P_{i}^{Gen})^{2}+\beta_{i}P_{i}^{Gen}+\gamma_{i},\qquad i\epsilon S_{G}\label{eq:cost function}
\end{equation}

A DG's surplus (commonly known as profit) measures the welfare of
a DG. In other words, it is a benefit gained by DG, when the sale
price of energy is more than the costs spent to produce the energy.
If the $i-th$ DG sells $\unit[P]{kW}$ of electrical energy at $\lambda\;\unitfrac{\$}{kWh}$,
the DG's profit is expressed as in (\ref{eq:DGsurplus}). 

\begin{equation}
\mathcal{S}_{i}^{Gen}\left(P_{i}^{Gen}\right)=\lambda P_{i}^{Gen}-C_{i}\left(P_{i}^{Gen}\right),\qquad i\epsilon S_{G}\label{eq:DGsurplus}
\end{equation}

\begin{figure}[tbh] \centering
\includegraphics[width=0.7\columnwidth]{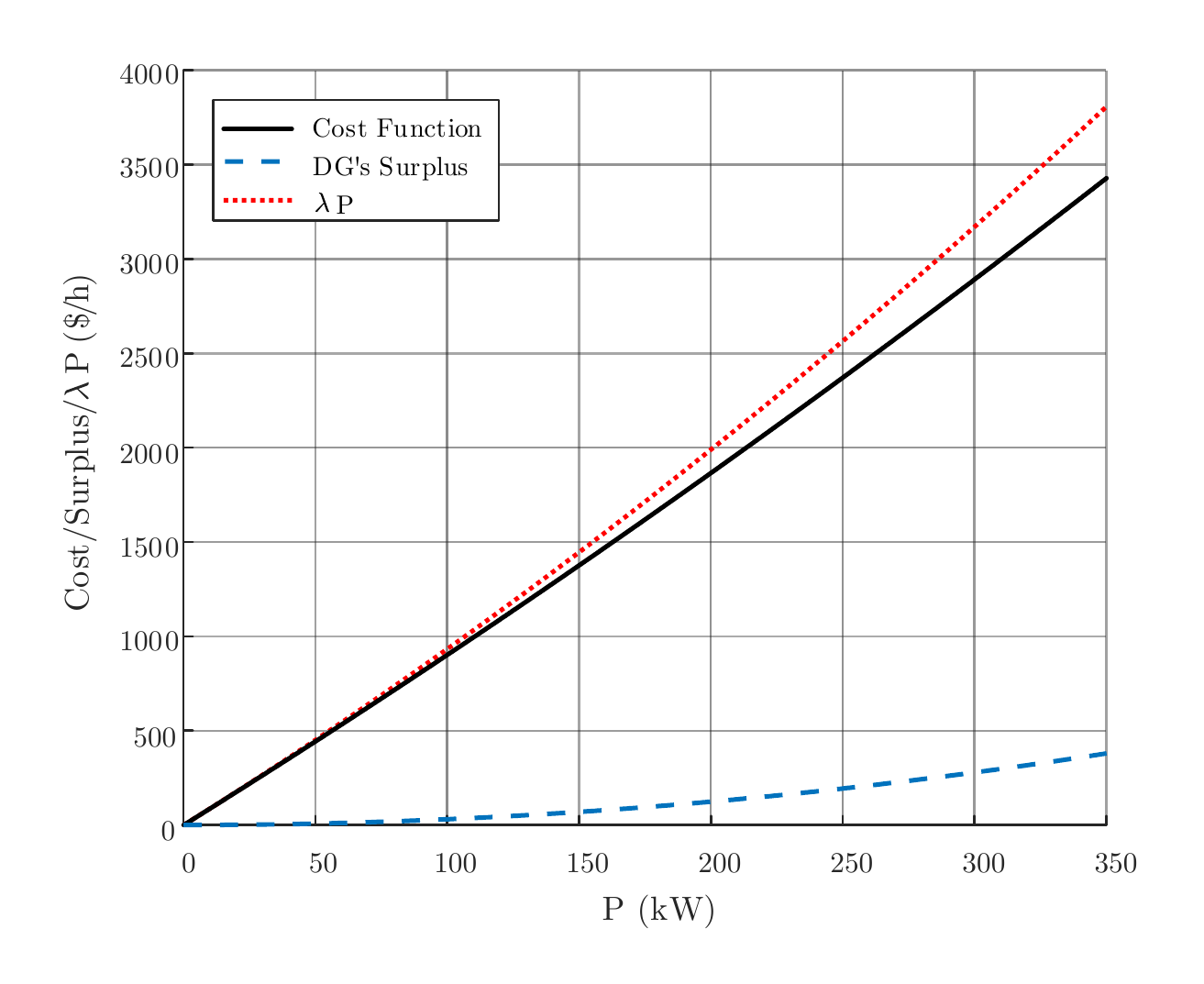}

\caption{A DG's surplus and cost function\label{fig:DGrsurpluscostfunction}}
\end{figure}

Unlike consumers, DGs tend to produce as much electricity as possible.
Figure \ref{fig:DGrsurpluscostfunction} shows the surplus curve of
a DG. The more power DGs produce, the higher surplus they obtain.

\subsection{Global Optimization Problem }

The objective function in this paper is to maximize the welfare of
all consumers and DGs. In other words, the objective function is to
maximize the summation of the utility functions (\ref{eq:utilityfunction})
and minimize the summation of the cost functions (\ref{eq:cost function}).
Thus, the overall objective function can be written as: 

\begin{gather}
Min\left(\underset{i\epsilon S_{G}}{\sum}C_{i}\left(P_{i}^{Gen}\right)-\underset{j\epsilon S_{D}}{\sum}U_{j}\left(P_{j}^{Load}\right)\right)\label{eq:OverallObjectiveFunction}
\end{gather}

Note that the current version of the proposed distributed algorithm
does not consider the power losses and the maximum capacity of the
power lines. This objective function is also subject to constraints
of power balance between the aggregated generations ($\underset{i\epsilon S_{G}}{\sum}P_{i}^{Gen}$
) and consumption ($\underset{j\epsilon S_{D}}{\sum}P_{j}^{Load}$
) as in (\ref{eq:powerbalance}).

\begin{equation}
\underset{i\epsilon S_{G}}{\sum}P_{i}^{Gen}=\underset{j\epsilon S_{D}}{\sum}P_{j}^{Load}\label{eq:powerbalance}
\end{equation}

In addition, the output-power of each DG and the consumption of each
consumer cannot go beyond their maximum capacity. These two constraints
can be applied to the optimization problem as in (\ref{eq:maximumcapacity}).

\begin{gather}
0\leq P_{i}^{Gen}\leq P_{i,max}^{Gen}\qquad i\epsilon S_{G}\nonumber \\
0\leq P_{j}^{Load}\leq P_{j,max}^{Load}\qquad j\epsilon S_{D}\label{eq:maximumcapacity}
\end{gather}

\section{Distributed Algorithm for Economic Dispatch\label{sec:DistributedAlgorithmforED}}

As previously mentioned, consensus-based distributed approaches offer
a great solution for solving optimization problems such as economic
dispatch. In this section, we review the conception of the graph theory
and elaborate on the proposed distributed algorithm in details.

\subsection{Graph Theory\label{subsec:GraphTheory}}

We can model the agents' interaction through the communication network
by (un)directed graphs denoted by $\mathcal{G}(\mathcal{V},\xi)$.
Consider a network of $n$ connected agents in which nodes are designated
by $\mathcal{V}=\{v_{1},v_{2},\ldots v_{p}\}$ and $\xi\subseteq\mathcal{V}\times\mathcal{V}$
shows a set of edge. The directed edge $\vec{e_{ij}}=(v_{i},v_{j})$
shows that agent $i$ share it's information state with agent $j$.
Also, undirected edge $e_{ij}=(v_{i},v_{j})$ indicates that agents
$i$ and $j$ can share information with each other. Two matrices
will commonly be used to represent the communication topology of a
multiple-agents network. The adjacency matrix denoted by $\mathcal{A}=\{\left[a_{ij}\right]|a_{ij}\in\mathcal{R^{P\times P}}\}$
of an undirected graph $\mathcal{G}$ is symmetric. The entry $a_{ij}$
of an adjacency matrix is a positive value if $e_{ij}\in\xi$ and
$a_{ij}=0$ for $e_{ij}\notin\xi$. Otherwise, the entry $a_{ii}$
is assumed to be zero. The second matrix is Laplacian matrix $L=\{\left[l_{ij}\right]|l_{ij}\in\mathcal{R^{P\times P}}\}$
in which entry $l_{ii}=\underset{j\quad}{\sum a_{ij}}$ and $l_{ij}=-a_{ij}$
for $i\neq j$ . Equation (\ref{eq:LandA}) shows matrices $\mathcal{A}$
and $L$ \cite{Soediono1989,Guo2014,Saber2003}.

\begin{multline}
A=\left[\begin{array}{cccc}
a_{11} & a_{12} & \cdots & a_{1n}\\
a_{21} & a_{22} & \cdots & a_{2n}\\
\vdots & \vdots & \ddots & \vdots\\
a_{n1} & a_{n2} & \cdots & a_{nn}
\end{array}\right]\\
L=\left[\begin{array}{cccc}
\underset{j\quad}{\sum a_{1j}} & -a_{12} & \cdots & -a_{1n}\\
-a_{21} & \underset{j\quad}{\sum a_{2j}} & \cdots & -a_{2n}\\
\vdots & \vdots & \ddots & \vdots\\
-a_{n1} & -a_{n2} & \cdots & \underset{j\quad}{\sum a_{nj}}
\end{array}\right]\qquad\qquad\quad\;\;\label{eq:LandA}
\end{multline}

\subsection{Consensus-based Distributed Protocols\label{subsec:ConsensusbasedDistributedProt}}

In consensus-based distributed approaches, a network of agents shares
information via communication channels between agents to reach a consensus.
Node $i$ and node $j$ have reached a consensus if and only if the
value of the state of the $i-th$ node ($x_{i}$) and the state of
the $j-th$ node ($x_{j}$) are equal \cite{Olfati-Saber2007,Saber2003}.
Thus, multiple agents reach a consensus when all of them agree on
the coordination information or variable. The Laplacian potential\textit{
}for a graph is delineated by (\ref{eq:laplacian}) which represents
a kind of virtual energy stored in a graph \cite{Melorose2015a}.

\begin{equation}
\mathcal{L_{P}=}\underset{i,j}{\sum}a_{ij}(x_{j}-x_{i})^{2}=2x^{T}Lx\label{eq:laplacian}
\end{equation}

In other words, the Laplacian potential could be used as a measure
that shows the total disagreement among all agents in a network. If
the agents of a network tend to reach a consensus, they should at
least interact with their neighbors to minimize Laplacian potential
($\mathcal{L_{P}}$) as a disagreement \cite{Saber2003}. In fact,
a general consensus for a multi-agent system will be reached if and
only if $\mathcal{L_{P}}=0$ or $x_{i}=x_{j}$.

The ``consensus'' being used for the proposed method is defined
as zero-power-mismatch. Based on the definition of Laplacian potential\textit{,
}the whole power mismatch is a virtual energy stored in the network
that must be minimized. Consensus is reached by converging towards

\begin{equation}
\varDelta P_{1}^{T}=\varDelta P_{2}^{T}=\cdots=\varDelta P_{n}^{T}=0
\end{equation}

where $\varDelta P_{i}^{T}$ is a power mismatch of the whole system
estimated by the $i-th$ DG.

Considering that all agents have single-integrator
dynamics \cite{Saber2003,Melorose2015a}, a standard linear consensus
protocol is defined as (\ref{eq:linearconsensus}).

\begin{equation}
\dot{x_{i}}(t)=\underset{j}{\sum}a_{ij}(x_{j}-x_{i})\label{eq:linearconsensus}
\end{equation}

Equation (\ref{eq:linearconsensus}) can be written for all agents as a vector: $\dot{x}(t)=-Lx$, which{}
is equivalent to the gradient of the Laplacian potential of a graph
as shown in Equation (\ref{eq:gradient}). It represents a gradient-descent
algorithm that is able to find the minimum of the Laplacian potential
function. As previously discussed, the minimum Laplacian
potential happens at $\varDelta P_{i}=0,\quad\;\forall i$.

\begin{equation}
\dot{x_{i}}(t)=-\nabla\mathcal{L_{P}}\label{eq:gradient}
\end{equation}

A discrete-time version of the linear consensus protocol of a first-order
integrator can be represented by (\ref{eq:discretetime})

\begin{gather}
x_{i}(k+1)-x(k)=u_{i}(k)\Longrightarrow\nonumber \\
x_{i}(k+1)=x_{i}(k)+u_{i}(k)\quad\quad\label{eq:discretetime}
\end{gather}

where $u_{i}(k)$ depends on the information state of neighbors of
$i-th$ agent (\textit{i.e.}, $u_{i}=f_{i}(x_{1},x_{2},...,x_{j})$)
and can be shown by (\ref{eq:uk}).

\begin{equation}
u_{i}(k)=\underset{j}{\sum}a_{ij}(x_{j}(k)-x_{i}(k))\label{eq:uk}
\end{equation}

Equation (\ref{eq:equation17}) is obtained from (\ref{eq:discretetime})
and (\ref{eq:uk}), where $I$ is unit matrix.

\begin{gather}
x_{i}(k+1)=x_{i}(k)+\underset{j}{\sum}a_{ij}(x_{j}(k)-x_{i}(k))\Longrightarrow\nonumber \\
x_{i}(k+1)=x_{i}(k)-Lx_{i}(k)=(I-L)x_{i}(k)\label{eq:equation17}
\end{gather}

Given that the sum of row of adjacency matrix A is one (\textit{i.e.},
A is a stochastic matrix), Equation (\ref{eq:finalprotocols}) can
be derived from (\ref{eq:equation17}). Equation (\ref{eq:finalprotocols})
explicitly indicates that the next state of each agent depends on
the current states of other agents.

\begin{align}
x_{i}(k+1) & =\underset{j}{\sum}a_{ij}x_{j}(k)\label{eq:finalprotocols}
\end{align}

Now, we consider $\varDelta P$ as information coordination that needs
to be shared among agents. Equation (\ref{eq:powerimbalanceasinfostate})
shows that $\varDelta P$ of each agent is calculated by the current
estimated $\varDelta P$ of other neighbors. It is worth mentioning
that $\varDelta P$, the only information shared among different agents,
does not include any private information. The consumers do not need to launch any communication link among themselves
for information coordination. Thus, the elements associated with the
connectivity between any pair of consumers are zero in \textquotedblleft Matrix
A\textquotedblright . Moreover, consumers do not have to establish
a communication channel with more than one DG. They are connected
to a local or nearest DG if there is a physical connection (transmission
and distribution line). Since the power mismatch is the only shared
information among DGs, a DG and its associated loads/consumers can
be viewed as an aggregate node. The \textquotedblleft Reduced Matrix
A\textquotedblright{} being used in Equation (19) is only an adjacency
matrix of DGs\textquoteright{} communication network. In sum, the
A matrix has many zero elements that most of them could be ignored
for the sake of simplicity. The A matrix that used in Equation (19)
is only an adjacency matrix of DGs\textquoteright{} communication
network without zero elements of consumers' network and communication
channel among consumers and DGs. Thus, we omitted zero elements. The
A matrix is reduced by dimension in comparison with A matrix of the
entire system.

\begin{multline}
\left[\begin{array}{c}
\varDelta P_{1}\left(k+1\right)\\
\varDelta P_{2}\left(k+1\right)\\
.\\
.\\
.\\
\varDelta P_{n}\left(k+1\right)
\end{array}\right]=\left[A\right]\left[\begin{array}{c}
\varDelta P_{1}\left(k\right)\\
\varDelta P_{2}\left(k\right)\\
.\\
.\\
.\\
\varDelta P_{n}\left(k\right)
\end{array}\right]\\
\\
\Longrightarrow\vec{\varDelta P}\left(k+1\right)=A\vec{\varDelta P}\left(k\right)\quad\quad\quad\quad\quad\quad\quad\label{eq:powerimbalanceasinfostate}
\end{multline}

Each DG uses its own estimated power mismatch as a feedback. By adding
the vector $\vec{P}^{L}(k)-\vec{P}^{Gen}\left(k\right)$ to
Equation (\ref{eq:powerimbalanceasinfostate}), Equation (\ref{eq:completedeltaP})
is obtained as a consensus protocol for this paper.

\begin{gather}
\vec{\varDelta P^{T}}\left(k+1\right)=A\vec{\varDelta P}\left(k\right)+\vec{P}^{L}\left(k\right)-\vec{P}^{Gen}\left(k\right)\nonumber \\
\nonumber \\
\Longrightarrow\left[\begin{array}{c}
\varDelta P_{1}^{T}\left(k+1\right)\\
\varDelta P_{2}^{T}\left(k+1\right)\\
.\\
.\\
.\\
\varDelta P_{n}^{T}\left(k+1\right)
\end{array}\right]=\left[A\right]\left[\begin{array}{c}
\varDelta P_{1}\left(k\right)\\
\varDelta P_{2}\left(k\right)\\
.\\
.\\
.\\
\varDelta P_{n}\left(k\right)
\end{array}\right]+\nonumber \\
\nonumber \\
\left[\begin{array}{c}
P_{1}^{L}\left(k\right)-P_{1}^{Gen}\left(k\right)\\
P_{2}^{L}\left(k\right)-P_{2}^{Gen}\left(k\right)\\
.\\
.\\
.\\
P_{n}^{L}\left(k\right)-P_{n}^{Gen}\left(k\right)
\end{array}\right]\qquad\nonumber \\
\label{eq:completedeltaP}
\end{gather}

Where $P_{i}^{L}$ is the summation of all local loads connected to $i-th$ DGs. Every iteration,
the DGs need to go through a simple process to update incremental
cost $\lambda$ internally.
This $\lambda$ does not need to be shared with neighbors. Once gain,
the only information that will be shared through the communication
network is $\varDelta P$.

The discrete-time equation (\ref{eq:DCformoflambda}) shows proposed
protocol for $\lambda$ in this paper, where $\varDelta x>0$ shows
the interval of discrete-time integration, and $K_{I}$ is the controller
coefficient.

\begin{equation}
\lambda_{i}\left(k+1\right)=K_{I}\underset{k}{\sum}\left(\varDelta P_{i}^{T}\left(k+1\right)+\varDelta P_{i}^{T}\left(k\right)\right)\times\frac{\varDelta x}{2}\label{eq:DCformoflambda}
\end{equation}

The incremental cost $\lambda$ is used to calculate output-power
($P_{i}^{Gen}$) of a DG, The parameters of cost function (\ref{eq:cost function})
such as $\alpha$, $\beta$ and the calculated $\lambda$ inside the
controller of each agent are used to estimate output-power using Equation
(\ref{eq:outputpower}). In fact, other agents and third-parties are
not at all able to access these parameters.

\begin{gather}
P_{i}^{Gen}\left(k+1\right)=\begin{cases}
0 & P_{i}^{Gen}\leq0\\
\frac{\lambda_{i}\left(k+1\right)-\beta_{i}}{2\alpha_{i}} & 0<P_{i}^{Gen}<P_{i,max}^{Gen}\\
P_{i,max}^{Gen} & P_{i}^{Gen}\geq P_{i,max}^{Gen}
\end{cases},\nonumber \\
\qquad i\epsilon S_{G}\label{eq:outputpower}
\end{gather}

When a DG estimates its output-power, it can determine the estimated
power mismatch by Equation (\ref{eq:powerimbalance}), and share this
estimate with its neighbors at each iteration.

\begin{equation}
\varDelta P_{i}\left(k+1\right)=P_{i}^{Gen}\left(k+1\right)-P_{j}^{L}(k+1),\qquad i\epsilon S_{G}\label{eq:powerimbalance}
\end{equation}

When a DG determines its $\lambda$ in accordance with its output-power,
it shares $\lambda$ with the local consumers. Then the consumers
calculate their demands based on the $\lambda$ offered by the DG.
It is not, however, necessary to share $\lambda$ among consumers;
in other words, each consumer just needs to receive $\lambda$ from
one DG. Then, the consumer can determine its best and most cost-effective
demand based on the maximum level of the consumer's surplus function
represented by (\ref{eq:consumersurplus}). The maximum of consumer's
surplus can be achieved by $\frac{\partial\mathcal{S}_{j}^{Load}}{\partial P_{j}^{Load}}=0$,
which is shown in (\ref{eq:maximumofconsumersurplus}). The utility
function of consumers in Equation (\ref{eq:utilityfunction}) and
Figure \ref{fig:Consumrsuplusutility} show that if a consumer uses
power more than $\nicefrac{\omega_{j}}{2b_{j}}$ , its level of satisfaction
will not be increased. Thus, the maximum load of $j-th$ consumer
is considered as $P_{j,max}^{Load}=\nicefrac{\omega_{j}}{2b_{j}}$.

\begin{gather}
P_{j}^{Load}(k+1)=\begin{cases}
0 & P_{j}^{Load}\leq0\\
\frac{\omega_{j}-\lambda_{i}\left(k+1\right)}{2b_{j}} & 0<P_{j}^{Load}<P_{j,max}^{Load}\\
P_{j,max}^{Load} & P_{j}^{Load}\geq P_{j,max}^{Load}
\end{cases},\nonumber \\
\qquad j\epsilon S_{D}\label{eq:maximumofconsumersurplus}
\end{gather}

The consumer sends the amount of estimated demand to a local DG if
there is a physical connection (i.e., distribution lines) between
them. The consumers do not need to disclose any properties of its
utility function. In addition, consumers do not need to establish
any communication channel among themselves to coordinate any information
state or with more than one DG. The above-mentioned features significantly
reduce the computing complexity and the upfront cost of new communication
infrastructure

\begin{figure}[tbh] \centering
\includegraphics[width=0.7\columnwidth]{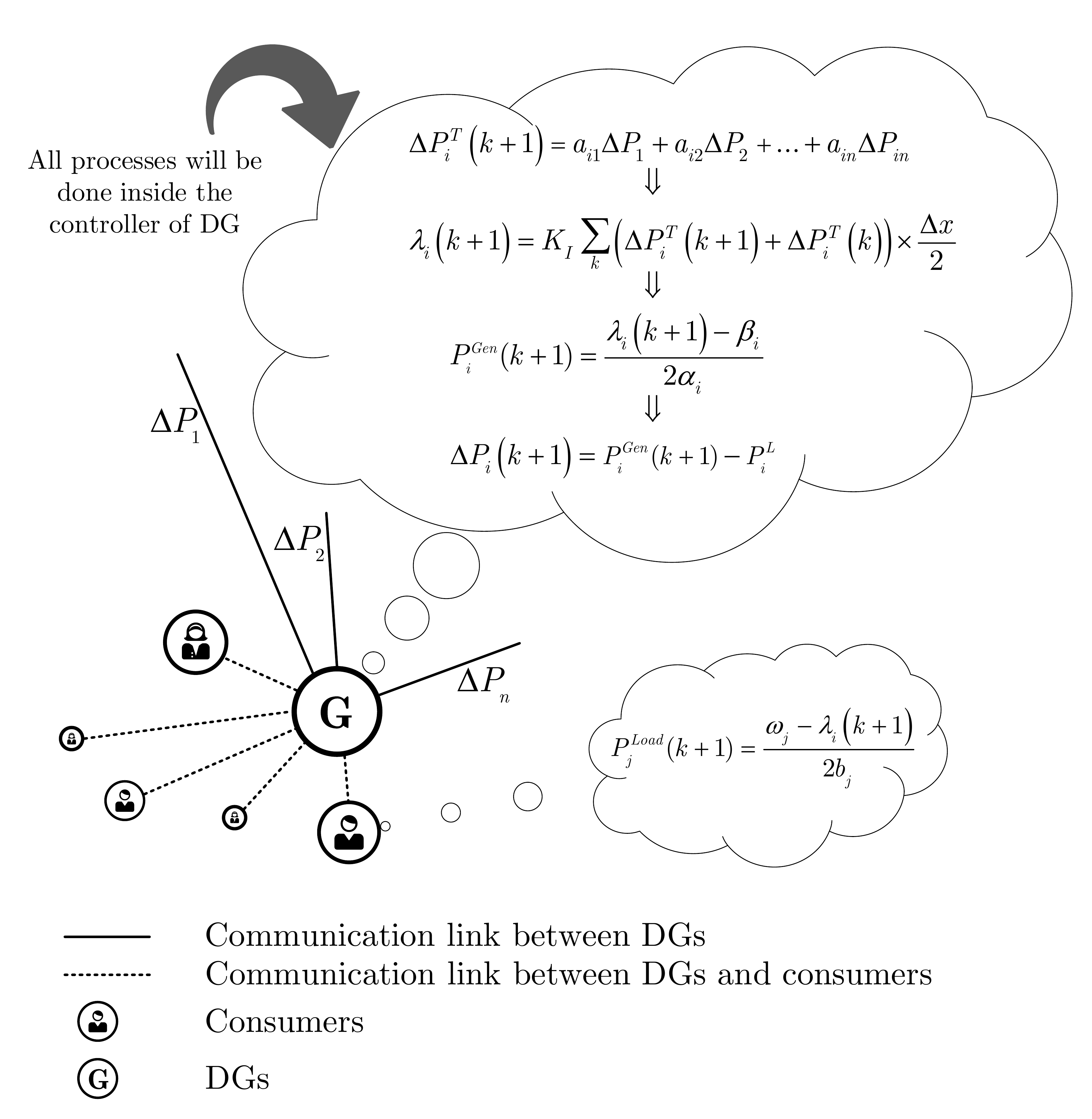}

\caption{Distributed decision making for a DG at every iteration.\label{fig:Interactionamongagents}}
\end{figure}

Figure \ref{fig:Interactionamongagents} illustrates the interaction
between a specific DG and other agents (DGs and consumers).

The optimality and convergence analysis of the proposed distributed
algorithm will be discussed in the rest of this section.

\subsection{The Optimality Analysis\label{subsec:OptimalityAnalysis}}

As mentioned before, $P_{i}^{L}$ is the summation of local loads
connected to $i-th$ DG. Equation (\ref{eq:summationoflocalloads})
shows the $\left(k+1\right)-th$ iteration of total $P_{i}^{L}$ and
$\psi_{i}$ indicates the portion of total load connected to $i-th$
DG. $\theta_{i}^{L}$ and $\varphi_{i}^{L}$ are used in the place
of $\psi_{i}\underset{j}{\sum}\frac{\omega_{j}}{2b_{j}}$ and $\psi_{i}\underset{j}{\sum}\frac{1}{2b_{j}}$
, respectively, for more simplicity. $\left(k+1\right)-th$ iteration
of DG's output is calculated by (\ref{eq:DGoutput}). In addition,
$\theta_{i}^{Gen}$ and $\varphi_{i}^{Gen}$ are used in the place
of $\frac{\beta_{i}}{2\alpha_{i}}$ and $\frac{1}{2\alpha_{i}}$,
respectively.

\begin{gather}
P_{i}^{L}(k+1)=\psi_{i}\underset{j}{\sum}P_{j}^{Load}(k+1)\nonumber \\
=\psi_{i}\underset{j}{\sum}\frac{\omega_{j}}{2b_{j}}-\psi_{i}\underset{j}{\sum}\frac{1}{2b_{j}}*\lambda_{i}(k+1)\nonumber \\
=\theta_{i}^{L}-\varphi_{i}^{L}\lambda_{i}(k+1);\quad\underset{j}{\sum}\psi_{i}=1\label{eq:summationoflocalloads}
\end{gather}

\begin{gather}
P_{i}^{Gen}(k+1)=\frac{\lambda_{i}(k+1)-\beta_{i}}{2\alpha_{i}}=\frac{\lambda_{i}(k+1)}{2\alpha_{i}}-\frac{\beta_{i}}{2\alpha_{i}}\nonumber \\
=\varphi_{i}^{Gen}\lambda_{i}(k+1)-\theta_{i}^{Gen}\label{eq:DGoutput}
\end{gather}

Finally, the power mismatch $\Delta P_{i}(k+1)$ calculated by $i-th$
DG is achieved by (\ref{eq:powermismatchcalculatedbyithDG}) where,
$\varphi_{i}^{c}=\varphi_{i}^{Gen}+\varphi_{i}^{L}$ and $\theta_{i}^{c}=\theta_{i}^{Gen}+\theta_{i}^{L}$.

\begin{gather}
\Delta P_{i}(k+1)=P_{i}^{Gen}(k+1)-P_{i}^{L}(k+1)\nonumber \\
=\left(\varphi_{i}^{Gen}+\varphi_{i}^{L}\right)\lambda_{i}(k+1)-\left(\theta_{i}^{Gen}+\theta_{i}^{L}\right)\nonumber \\
=\varphi_{i}^{c}\lambda_{i}(k+1)-\theta_{i}^{c}\label{eq:powermismatchcalculatedbyithDG}
\end{gather}

$\lambda_{i}$, $\Delta P_{i}^{T}$ of $i-th$ agent for $k=0$ can
be calculated by (20), (23) and (25) .

\begin{gather}
\Delta P_{i}^{T}(1)=\Delta P_{1}(0)+\Delta P_{2}(0)+\cdots=\Delta_{1}\\
\Delta P_{i}^{T}(0)=\Delta P_{1}(0^{-})+\Delta P_{2}(0^{-})+\cdots=\Delta_{0}\\
\lambda_{i}(1)=K_{I}(\Delta P_{i}^{T}(1)+\Delta P_{i}^{T}(0))=K_{I}(\Delta_{1}+\Delta_{0})
\end{gather}

By continuing this process for $k=1,\,2,\,\cdots$, equations (\ref{eq:FinaldeltaP}),
(\ref{eq:FinalLambda}) will be obtained.

\begin{equation}
\Delta_{k+1}=\Delta P_{i}^{T}(k+1)=\underset{i}{\sum}K_{I}\varphi_{i}^{c}\Delta_{k}-\Delta_{k}\quad\forall i\epsilon S_{G}\label{eq:FinaldeltaP}
\end{equation}

\begin{equation}
\lambda_{i}(k+1)=K_{I}(\Delta_{k}+\cdots+2\Delta_{2}+2\Delta_{1}+\Delta_{0})\quad\forall i\epsilon S_{G}\label{eq:FinalLambda}
\end{equation}

In Equation (\ref{eq:FinaldeltaP}), $||\Delta_{k+1}||$ merges to
zero \textit{i.e.; }$||\Delta_{k+1}||\rightarrow0$ if $||\underset{i}{1-\sum}K_{I}\varphi_{i}^{c}||\leq\varepsilon<1$,
where where $\varepsilon$ is a positive number. Then, we have $\lambda_{i}(k+1)\leq(\varepsilon^{k}\Delta_{0}+2\varepsilon^{k-1}\Delta_{0}\cdots+2\varepsilon^{2}\Delta_{0}+2\varepsilon\Delta_{0}+\Delta_{0}),\quad\forall i\epsilon S_{G}$.
Therefore, if $k\rightarrow\infty$, $\lambda_{i}$will not approach
infinity; hence, $\lambda_{i}(k+1)\leq\frac{1}{1-\varepsilon}\Delta_{0}+\frac{\varepsilon}{1-\varepsilon}\Delta_{0}=\frac{1+\varepsilon}{1-\varepsilon}\Delta_{0}$.
The smaller $||\underset{i}{1-\sum}K_{I}\varphi_{i}^{c}||$ is, the
faster the power mismatch will converge to $0$. $K_{I}$ is the parameter
that can control/change the size of $||\underset{i}{1-\sum}K_{I}\varphi_{i}^{c}||$
to be less than $1$.

Finally, Equation (\ref{eq:finalconvergedsolution}) indicates that
$\lambda_{i}$ for all DG ($\forall i$) will be same and $\Delta P_{i}^{T}(k+1)$
of $i-th$ of DG will converge to zero for $\forall i\epsilon S_{G}$.
We consider all $\lambda_{i}\quad\forall i$ as $\lambda$ because
they are identical. In the next step, we will show (\ref{eq:finalconvergedsolution})
will satisfy the KKT conditions.

\begin{gather}
\underset{k\rightarrow\infty}{lim}\;\lambda_{i}(k+1)=\lambda;\qquad\forall i\epsilon S_{G}\nonumber \\
\underset{k\rightarrow\infty}{lim}\;\Delta P_{i}^{T}(k+1)=\Delta_{k+1}=0;\qquad\forall i\epsilon S_{G}\label{eq:finalconvergedsolution}
\end{gather}

\textbf{Assumption 1}: All local cost functions utility functions
(\ref{eq:utilityfunction}) and (\ref{eq:cost function}) are strictly
concave and convex, respectively. Accordingly, the total objective
function (\ref{eq:OverallObjectiveFunction}) is strictly convex.

\textbf{Assumption 2}: In addition, all equality and inequality constraint
functions, represented by (\ref{eq:powerbalance}) and (\ref{eq:maximumcapacity}),
are affine.

\textbf{Lemma 1}: The optimization problem represented in this paper
through (\ref{eq:utilityfunction})-(\ref{eq:maximumcapacity}) is
a convex optimization problem with differentiable objective and constraint
functions satisfying Slater\textquoteright s condition (\textit{assumption
1 and 2 guarantee Slater's condition}), thus the KKT conditions provide\textbf{\textit{
necessary and sufficient conditions}} for optimality \cite{Boyd2004}.

The remaining of this section makes certain that the fixed-point of
proposed iterative consensus algorithm obtained by (\ref{eq:completedeltaP})-(\ref{eq:maximumofconsumersurplus})
is a global optimal solution of (\ref{eq:OverallObjectiveFunction})
if it is satisfying the following KKT conditions \cite{Yi2015}.

\noindent \textbf{Lagrangian:}

\begin{gather}
L(P,\lambda,\mu,\zeta)=\left(\underset{i\epsilon S_{G}}{\sum}C_{i}\left(P_{i}^{Gen}\right)-\underset{j\epsilon S_{D}}{\sum}U_{j}\left(P_{j}^{Load}\right)\right)\nonumber \\
+\lambda\left(\underset{j\epsilon S_{D}}{\sum}P_{j}^{Load}-\underset{i\epsilon S_{G}}{\sum}P_{i}^{Gen}\right)\nonumber \\
+\underset{i\epsilon S_{G}}{\sum}\mu_{i}\left(P_{i}^{Gen}-P_{i,max}^{Gen}\right)+\underset{i\epsilon S_{G}}{\sum}\zeta_{i}\left(-P_{i}^{Gen}\right)\nonumber \\
+\underset{j\epsilon S_{D}}{\sum}\mu_{j}\left(P_{j}^{Load}-P_{j,max}^{Load}\right)+\underset{j\epsilon S_{D}}{\sum}\zeta_{j}\left(-P_{j}^{Load}\right)
\end{gather}

\noindent \textbf{Lagrangian stationarity }($\nabla_{P}L(P,\lambda,\mu,\zeta)$=0):

\begin{gather}
\frac{\partial C_{i}\left(P_{i}^{Gen}\right)}{\partial\left(P_{i}^{Gen}\right)}-\lambda+\mu_{i}-\zeta_{i}=0;\qquad\forall i\epsilon S_{G}\nonumber \\
-\frac{\partial U_{j}\left(P_{j}^{Load}\right)}{\partial\left(P_{j}^{Load}\right)}+\lambda+\mu_{j}-\zeta_{j}=0;\qquad\forall j\epsilon S_{D}\label{eq:Lagrangianstationarity}
\end{gather}

\noindent \textbf{Complementary slackness:}

\begin{gather}
\mu_{i}\left(P_{i}^{Gen}-P_{i,max}^{Gen}\right)=0\,,\quad\zeta_{i}\left(-P_{i}^{Gen}\right)=0;\qquad\forall i\epsilon S_{G}\nonumber \\
\mu_{j}\left(P_{j}^{Load}-P_{j,max}^{Load}\right)=0\,,\quad\zeta_{j}\left(-P_{j}^{Load}\right)=0;\qquad\forall j\epsilon S_{D}\label{eq:Complementaryslackness}
\end{gather}

\noindent \textbf{Dual feasibility:}

\begin{gather}
\mu_{i}\geq0\,,\quad\zeta_{i}\geq0;\qquad\forall i\epsilon S_{G}\nonumber \\
\mu_{j}\geq0\,,\quad\zeta_{j}\geq0;\qquad\forall j\epsilon S_{D}\label{eq:Dualfeasibility}
\end{gather}

All local constraints presented by (\ref{eq:maximumcapacity}) for
generation and consumption level of DGs and loads are considered as
the primal feasibility.

Let consider fixed point of the proposed iterative consensus algorithm
as optimal point, $\Delta P_{i}^{T}=\Delta P_{1}+\Delta P_{2}+\cdots\Delta P_{n}=0$,
where $\Delta P_{i}=P_{i}^{Gen}-P_{i}^{L}$ and $\lambda_{1}=\lambda_{2}=\cdots=\lambda_{n}=\lambda,\;\forall i$.
$\Delta P_{i}^{T}=0$ satisfies the equality constraint, as the load balance of the entire
system, to ensure that the demand will be supported. As mentioned
before, there is only one equality constraint (\ref{eq:powerbalance});
thus, all agents should reach the same $\lambda$ and (\ref{eq:summationoflocalloads})-(\ref{eq:finalconvergedsolution})
guarantee the identical $\lambda$ for all agents.

If the optimal solution of the objective function does not violate
local constraints (inequality functions) represented by (\ref{eq:maximumcapacity}),
then these constraints will never play any role in changing the minimum
compared with the problem without the inequality constraints. The
DG's profit is maximized when $\nicefrac{\partial\mathcal{S}_{i}^{Gen}}{\partial P_{i}^{Gen}}=0;\:\forall i\epsilon S_{G}$.
It means $\lambda_{opt}-\nicefrac{\partial C_{i}}{\partial P_{i}^{Gen}}|_{P_{i,opt}^{Gen}}=0;\:\forall i\epsilon S_{G}$.
The Lagrangian stationarity (\ref{eq:Lagrangianstationarity}), complementary
slackness (\ref{eq:Complementaryslackness}) and dual feasibility
(\ref{eq:Dualfeasibility}) are satisfied by taking $\lambda=\nicefrac{\partial C_{i}}{\partial P_{i}^{Gen}}|_{P_{i,opt}^{Gen}};\:\forall i\epsilon S_{G}$ and $\mu_{i}=\zeta_{i}=0,\;\forall i$. Thus, $\lambda$ obtained
by algorithm is $\lambda_{opt}$ because it satisfies the KKT condition.

The local constraints can affect the optimal solution in two ways:
\begin{itemize}
\item If the optimal points are greater than maximum limitation, $P_{i,opt}^{Gen}\geq P_{i,max}^{Gen}$,
in this case, the maximum level is considered as an optimal solution,
hence $P_{i,opt}^{Gen}=P_{i,max}^{Gen}$ and $\nicefrac{\partial C_{i}}{\partial P_{i}^{Gen}}|_{P_{i,max}^{Gen}}-\lambda_{opt}\leq0$.
Thus, the Lagrangian stationarity (\ref{eq:Lagrangianstationarity}),
complementary slackness (\ref{eq:Complementaryslackness}) and dual
feasibility (\ref{eq:Dualfeasibility}) are satisfied by taking $\lambda=\lambda_{opt};\;\forall i,j$,
$\mu_{i}=\lambda_{opt}-\nicefrac{\partial C_{i}}{\partial P_{i}^{Gen}}|_{P_{i,max}^{Gen}}$
and $\zeta_{i}=0,\;\forall i$.
\item If the optimal points are less than minimum limitation, $P_{i,opt}^{Gen}\leq0$,
in this case, the minimum level is considered as optimal solution,
hence $P_{i,opt}^{Gen}=0$ and $\nicefrac{\partial C_{i}}{\partial P_{i}^{Gen}}|_{0}-\lambda_{opt}\geq0$.
Thus, the Lagrangian stationarity (\ref{eq:Lagrangianstationarity}),
complementary slackness (\ref{eq:Complementaryslackness}) and dual
feasibility (\ref{eq:Dualfeasibility}) are satisfied by taking $\lambda=\lambda_{opt};\;\forall i,j$,
$\zeta_{i}=\nicefrac{\partial C_{i}}{\partial P_{i}^{Gen}}|_{0}-\lambda_{opt}$
and $\mu_{i}=0,\;\forall i$.
\end{itemize}
The same procedure could be considered for consumers' side. If all
local constraints (inequality functions) represented by (9) are ignored,
the consumer surplus is maximized when $\nicefrac{\partial\mathcal{S}_{j}^{Load}}{\partial P_{j}^{Load}}=0;\:\forall j\epsilon S_{D}$.
It means $\nicefrac{\partial U_{j}}{\partial P_{j}^{Load}}|_{P_{j,opt}^{Load}}-\lambda_{opt}=0;\:\forall j\epsilon S_{D}$.
The Lagrangian stationarity (\ref{eq:Lagrangianstationarity}), complementary
slackness (\ref{eq:Complementaryslackness}) and dual feasibility
(\ref{eq:Dualfeasibility}) are satisfied by taking $\lambda=\nicefrac{\partial U_{j}}{\partial P_{j}^{Load}}|_{P_{j,opt}^{Load}};\:\forall j\epsilon S_{D}$
and $\mu_{j}=\zeta_{j}=0,\;\forall i$. Thus, $\lambda$ obtained
by algorithm is $\lambda_{opt}$ because it satisfies the KKT condition.
\begin{itemize}
\item If the optimal points are greater than maximum limitation, $P_{j,opt}^{Load}\geq P_{j,max}^{Load}$,
in this case, the maximum level is considered as an optimal solution,
hence $P_{j,opt}^{Load}=P_{j,max}^{Load}$ and $\nicefrac{\partial U_{j}}{\partial P_{j}^{Load}}|_{P_{i,max}^{Load}}-\lambda_{opt}\geq0$.
Thus, the Lagrangian stationarity (\ref{eq:Lagrangianstationarity}),
complementary slackness (\ref{eq:Complementaryslackness}) and dual
feasibility (\ref{eq:Dualfeasibility}) are satisfied by taking $\lambda=\lambda_{opt};\;\forall j$,
$\mu_{j}=\nicefrac{\partial U_{j}}{\partial P_{j}^{Load}}|_{P_{i,max}^{Load}}-\lambda_{opt}$
and $\zeta_{j}=0,\;\forall j$.
\item If the optimal points are less than minimum limitation, $P_{j,opt}^{Load}\leq0$,
in this case, the minimum level is considered as optimal solution,
hence $P_{j,opt}^{Load}=0$ and $\nicefrac{\partial U_{j}}{\partial P_{j}^{Load}}|_{0}-\lambda_{opt}\leq0$.
Thus, the Lagrangian stationarity (\ref{eq:Lagrangianstationarity}),
complementary slackness (\ref{eq:Complementaryslackness}) and dual
feasibility (\ref{eq:Dualfeasibility}) are satisfied by taking $\lambda=\lambda_{opt};\;\forall i,j$,
$\zeta_{j}=\lambda_{opt}-\nicefrac{\partial U_{j}}{\partial P_{j}^{Load}}|_{0}$
and $\mu_{j}=0,\;\forall j$.
\end{itemize}
\textbf{Lemma 2}: The optimal solution found by the proposed iterative
consensus algorithm is unique.

It is proved that the fixed point of the proposed iterative consensus
algorithm satisfies the KKT conditions while objective and constraint
functions are both strictly convex and differentiable. Thus, the satisfaction
of Slater's condition provides an absolute assurance about global
optimality \cite{Ruszczynski2011}.

\section{Performance Assessment\label{sec:PerformanceAssessment}}

In this section, we conduct a performance evaluation of the proposed
distributed method through three case studies using software simulations
and experimental test. All software simulations are conducted in the
MATLAB 2015a environment on an ordinary desktop PC with an Intel(R)
Core(TM)i3 CPU @ 2.13 GHz, 4-GB RAM memory. The experiment test is
performed using the VOLTTRON\texttrademark{} platform
and a cluster of low-cost credit-card-size single-board PCs.

In the first case study, we provide a numerical example to evaluate
the algorithm performance (i.e., accuracy) in a relatively small-scale
system. The numerical results are compared with the benchmark results
found by a traditional centralized economic dispatch. The centralized
method is implemented using YALMIP toolbox \cite{YALMIPTOOLBOX} and
MATLAB.

In the second case study, we demonstrate the scalability and fast
convergence rate of the proposed distributed algorithm in a large-scale
network with 1,400 agents.

\begin{figure}[tbh] \centering
\includegraphics[width=0.7\columnwidth]{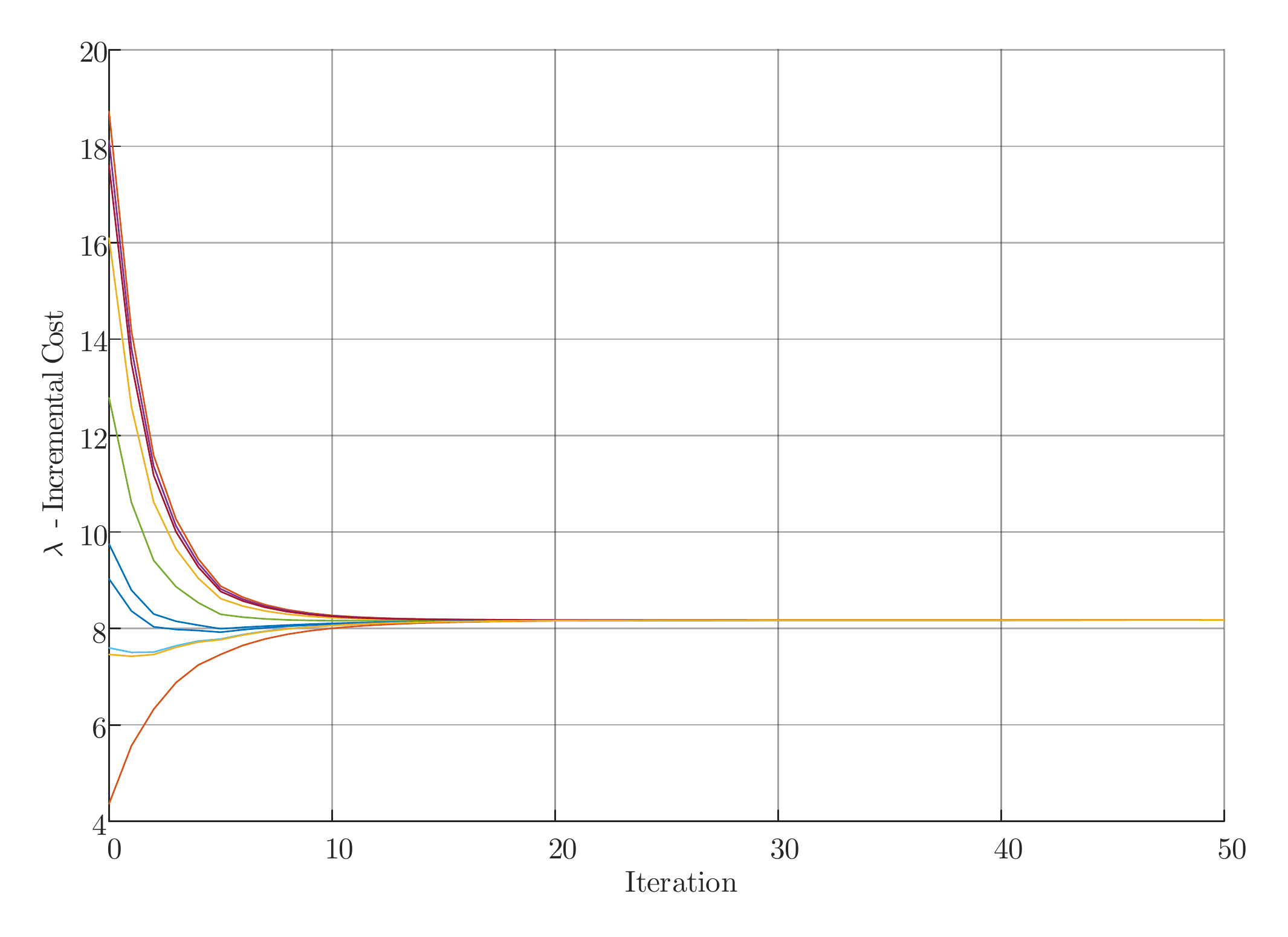}

\caption{Incremental cost of DGs in case study I \label{fig:Incrementalcostof}}
\end{figure}

\begin{figure}[tbh] \centering
\includegraphics[width=0.7\columnwidth]{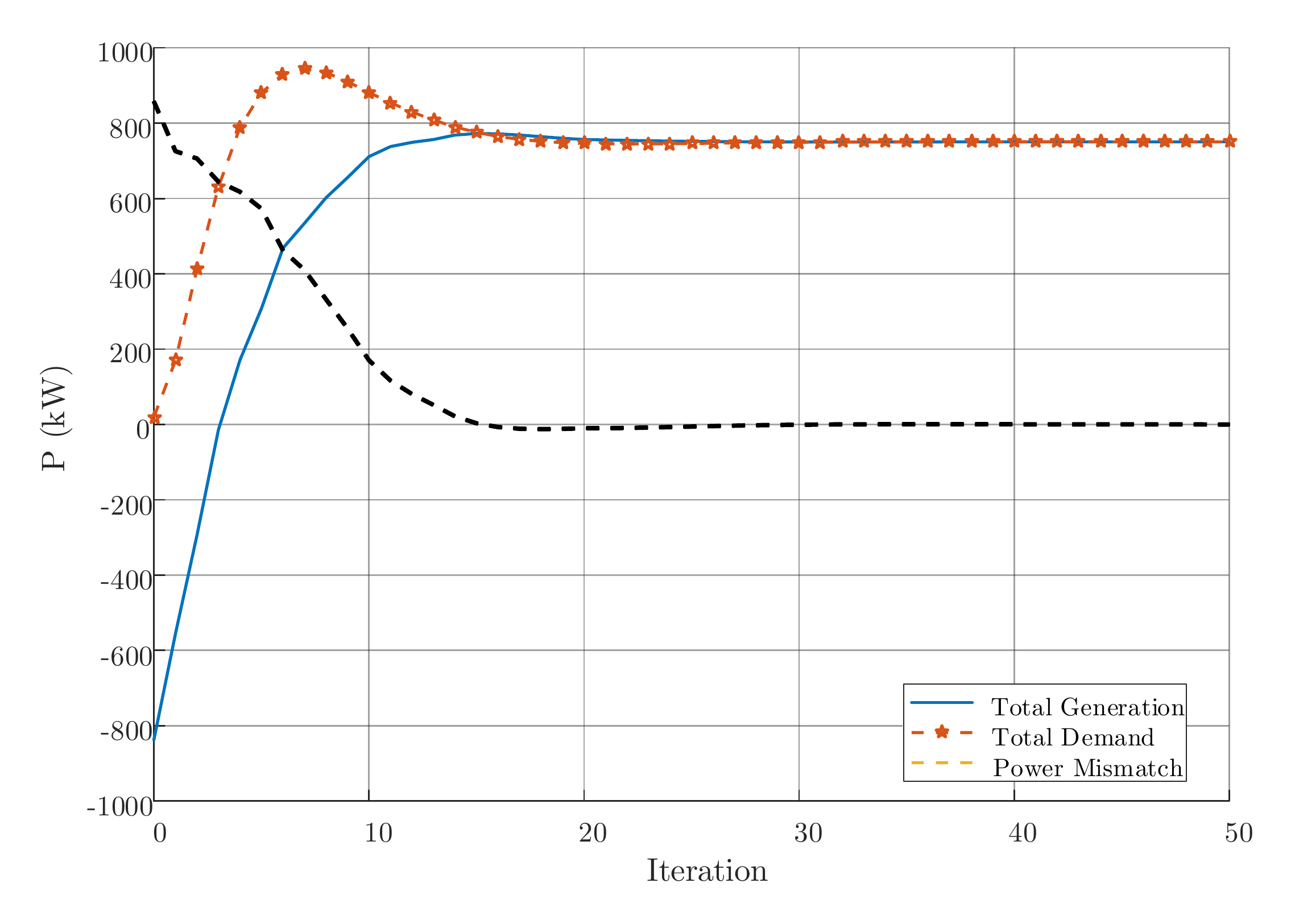}

\caption{Generation, demand and power mismatch in case\label{fig:Generationdemandandpowerimbalance}}
\end{figure}

In the third case study, we set up an experimental testbed to verify
the practical performance of the proposed distributed algorithm using
the {VOLTTRON\texttrademark{} platform and a group
of low-cost credit-card-size single-board PCs.

\subsection{Case Study I (Evaluation of Accuracy)\label{subsec:CaseStudyI}}

In this case study, an IEEE 39-bus test system with 29 agents (10
DGs and 19 consumers) is considered. Their cost functions and utility
functions are formulated using Equations (\ref{eq:cost function})
and (\ref{eq:utilityfunction}), respectively. Table \ref{tab:parametersofagents}
summarizes the parameters of the cost functions and utility functions
of the agents \cite{Rahbari-Asr2014a}. The initial values of the
$\lambda$s are randomly selected. The controller parameters ($K_{I}$
and $K_{P}$) are obtained by trial-and-error. $K_{I}$ and $K_{P}$
can be randomly set in the range of zero to one.

\begin{figure}[tbh] \centering
\includegraphics[width=0.7\columnwidth]{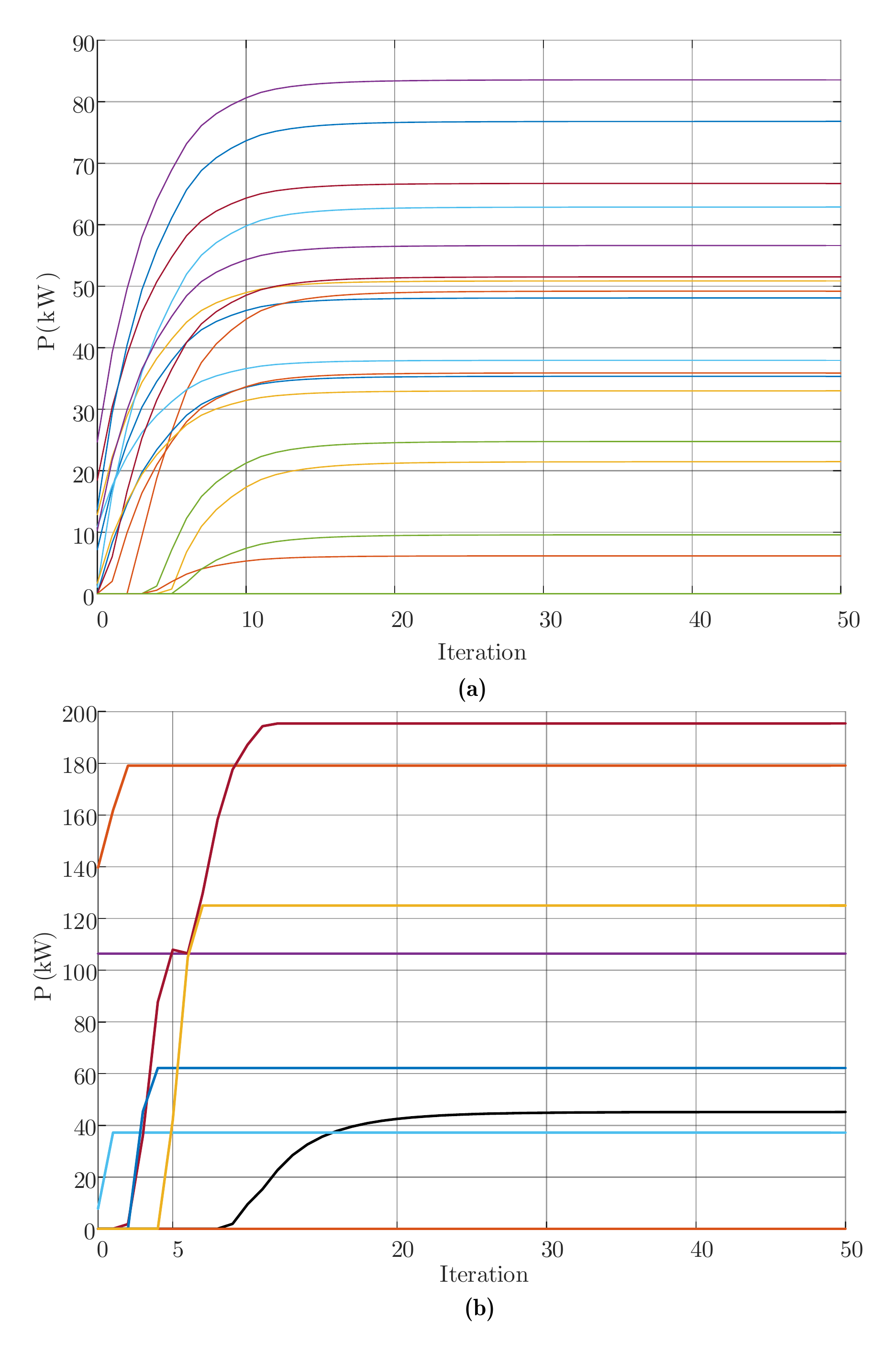}

\caption{Generation output and load demand powers in case study I; (a) demands
($kW$), (b) DGs ($kW$)\label{fig:GenerationandLoad}}
\end{figure}

\begin{table}[tbh]
\caption{The coefficients of DGs' cost functions and consumers' utility functions
in case study I\label{tab:parametersofagents}}

\center

\begin{tabular}{c|ccc|ccc}
\multicolumn{1}{c}{} & \multicolumn{3}{c}{Cost function} & \multicolumn{3}{c}{Utility function}\tabularnewline
\hline
$i$ & $\alpha_{i}$ & $\beta_{i}$ & \multicolumn{1}{c}{$P_{i,max}^{Gen}$} & $\omega_{i}$ & $b_{i}$ & $P_{i,max}^{Load}$\tabularnewline
\hline
1 & 0.0031  & 8.71  & 113.23 & 17.17  & 0.0935 & 91.79 \tabularnewline
2 & 0.0074  & 3.53  & 179.1  & 12.28  & 0.0417 & 147.29 \tabularnewline
3 & 0.0066  & 7.58  & 90.03  & 18.42 & 0.1007  & 91.41 \tabularnewline
4 & 0.0063  & 2.24  & 106.41  & 7.06  & 0.0561 & 62.96 \tabularnewline
5 & 0.0069  & 8.53  & 193.80  & 10.85  & 0.0540 & 100.53 \tabularnewline
6 & 0.0014  & 2.25 & 37.19  & 18.91  & 0.1414 & 66.88 \tabularnewline
7 & 0.0041  & 6.29  & 195.4 & 18.76  & 0.0793 & 118.35 \tabularnewline
8 & 0.0051 & 4.30  & 62.17  & 15.70  & 0.1064 & 73.81 \tabularnewline
9 & 0.0032  & 8.26  & 143.41  & 14.28  & 0.0850 & 84.00\tabularnewline
10 & 0.0025 & 5.3 & 125 & 10.15 & 0.0460  & 110.32 \tabularnewline
11 & \textendash{} & \textendash{} & \textendash{} & 19.04  & 0.0650  & 146.46 \tabularnewline
12 & \textendash{} & \textendash{} & \textendash{} & 06.87  & 0.0549  & 62.61 \tabularnewline
13 & \textendash{} & \textendash{} & \textendash{} & 15.96 & 0.0619 & 128.91 \tabularnewline
14 & \textendash{} & \textendash{} & \textendash{} & 14.70  & 0.0633  & 116.08 \tabularnewline
15 & \textendash{} & \textendash{} & \textendash{} & 17.50  & 0.0607  & 144.04 \tabularnewline
16 & \textendash{} & \textendash{} & \textendash{} & 10.97  & 0.2272 & 24.15 \tabularnewline
17 & \textendash{} & \textendash{} & \textendash{} & 16.25 & 0.1224 & 66.39 \tabularnewline
18 & \textendash{} & \textendash{} & \textendash{} & 17.53  & 0.0826  & 106.14 \tabularnewline
19 & \textendash{} & \textendash{} & \textendash{} & 09.84 & 0.0869 & 56.60\tabularnewline
\end{tabular}
\end{table}

\begin{table}[tbh]
\caption{Comparison of the output of DGs achieved by distributed and centralized
methods in case study I\label{tab:Comparisonofgenerators}}

\center

\begin{tabular}{c|cc}
\hline
Output Power $(kW)$ & Distributed Method & Centralized Method\tabularnewline
$P_{1}$ & 0 & 0\tabularnewline
$P_{2}$ & 179.1 & 179.099\tabularnewline
$P_{3}$ & 45.16 & 45.1614\tabularnewline
$P_{4}$ & 106.4 & 106.409\tabularnewline
$P_{5}$ & 0 & 0\tabularnewline
$P_{6}$ & 37.19 & 37.189\tabularnewline
$P_{7}$ & 195.4 & 195.399\tabularnewline
$P_{8}$ & 62.17 & 62.1699\tabularnewline
$P_{9}$ & 0 & 0\tabularnewline
$P_{10}$ & 125 & 124.999\tabularnewline
Total & 750.4 & 750.428\tabularnewline
\end{tabular}
\end{table}

Figure \ref{fig:GenerationandLoad} shows the evolution of DG power
output and load demand, respectively. Figure \ref{fig:GenerationandLoad}(a)
contains 19 consumer demand curves and Figure \ref{fig:GenerationandLoad}(b)
includes 10 DG output-power curves. As can be seen in \ref{fig:GenerationandLoad},
the economic dispatch solution converges at the 36-th iteration. The
corresponding execution time is about 1.69 seconds. The accuracy of
the proposed distributed solution algorithm is validated by the benchmark
results found by a centralized method. As shown in Tables
\ref{tab:Comparisonofgenerators} and \ref{tab:Comparisonofloads},
the solution mismatch between the distributed and centralized methods
is less than 0.00201\% of the average. As the distributed algorithm
proceeds, the incremental cost converges to $\unitfrac[8.175]{\$}{kwh}$,
as shown in Figure \ref{fig:Incrementalcostof}.
The evolution of power mismatch, the evolution of
total generation and the evolution of total load demand are shown
in Figure \ref{fig:Generationdemandandpowerimbalance}. Power mismatch
($\Delta P$) serves as \textit{coordination information }and
gradually converges to a consensus \textit{(i.e., $0\;kW$). }The
power tolerance is set to be 0.001 kW in our case studies.
As the incremental cost and power mismatch settle down, the optimal
value of social welfare is found to be $\unitfrac[5,211.5]{\$}{hr}$.
It is important to note that the proposed distributed
algorithm is able to converge to the near-optima much faster than
other distributed methods \cite{Hug2015a,Guo2014,Rahbari-Asr2014a}.
For example, one of the published works \cite{Rahbari-Asr2014a} showed
that the same economic dispatch problem was solved after 500 iterations,
while our distributed control algorithm is able to find the same results
at the 36-th iteration.

\begin{table}[tbh]
\caption{Comparison of the demand of loads achieved by distributed and centralized
methods\label{tab:Comparisonofloads}}

\center

\begin{tabular}{c|cc}
\hline
Demand of Loads $(kW)$ & Distributed Method & Centralized Method\tabularnewline
$L_{1}$ & 48.1 & 48.095\tabularnewline
$L_{2}$ & 49.21 & 49.207\tabularnewline
$L_{3}$ & 50.86 & 50.863\tabularnewline
$L_{4}$ & 0 & 0\tabularnewline
$L_{5}$ & 24.76 & 24.758\tabularnewline
$L_{6}$ & 37.96 & 37.955\tabularnewline
$P_{7}$ & 66.73 & 66.733\tabularnewline
$L_{8}$ & 35.36 & 35.356\tabularnewline
$L_{9}$ & 35.91 & 35.905\tabularnewline
$L_{10}$ & 21.46 & 21.455\tabularnewline
$L_{11}$ & 83.57 & 83.568\tabularnewline
$L_{12}$ & 0 & 0\tabularnewline
$L_{13}$ & 62.87 &  62.874\tabularnewline
$L_{14}$ & 51.53 & 51.531\tabularnewline
$L_{15}$ & 76.8 & 76.80\tabularnewline
$L_{16}$ & 6.148 & 6.148\tabularnewline
$L_{17}$ & 32.98 & 32.981\tabularnewline
$L_{18}$ & 56.62 & 56.621\tabularnewline
$L_{19}$ & 9.573 & 9.573\tabularnewline
Total & 750.4 & 750.428\tabularnewline
\end{tabular}
\end{table}

\subsection{Case Study II (Evaluation of Scalability and Fast Convergence)}

In order to demonstrate the scalability, we then apply the proposed
solution algorithm to a large-scale system, including
1,000 consumers and 400 DGs. The initial conditions of the1,400 agents
are randomly selected.

As shown in Figure \ref{fig:1400v2}, the incremental
costs reach consensus within approximately 40 iterations, which is
considered a fast convergence rate for a 1,400-agent system. The corresponding
execution time is 192.579 seconds. Some simulation for different numbers
of agents, 29 (10 DGs and 19 Consumers), 350 (150 DGs and 200 Consumers),
700 (300 DGs and 400 Consumers), 1050 (350 DGs and 700 Consumers),
1400 (400 DGs and 1000 Consumers), are performed to study the trend
in number of iteration for convergence. The Figure \ref{fig:Thetrend}
shows that as the number of agents increase from 29 (case study I)
to 1,400 (case study II), the number of iteration is almost constant,
demonstrating that the proposed distributed algorithm is particularly
suitable for solving large-scale economic dispatch problems. The minimum
error criteria for power mismatch tolerance is set to be}\textit{
$\Delta P=0.001\;kW$ }in our case studies used in Figure
\ref{fig:Thetrend}.

As previously emphasized, the power mismatch is the only shared information
between agents. It reduces the computational cost because the proposed
algorithm only needs a simple updating process on the power mismatch.
Besides, consumers do not need to establish any communication channel
among themselves to coordinate any information state or with more
than one DG. The above-mentioned features contribute to a significantly
reduction on the computing complexity.

\subsection{Case Study III (Evaluation of Practical Performance)\label{subsec:CaseStudyIII}}

\begin{figure}[tbh] \centering
\includegraphics[width=0.5\columnwidth]{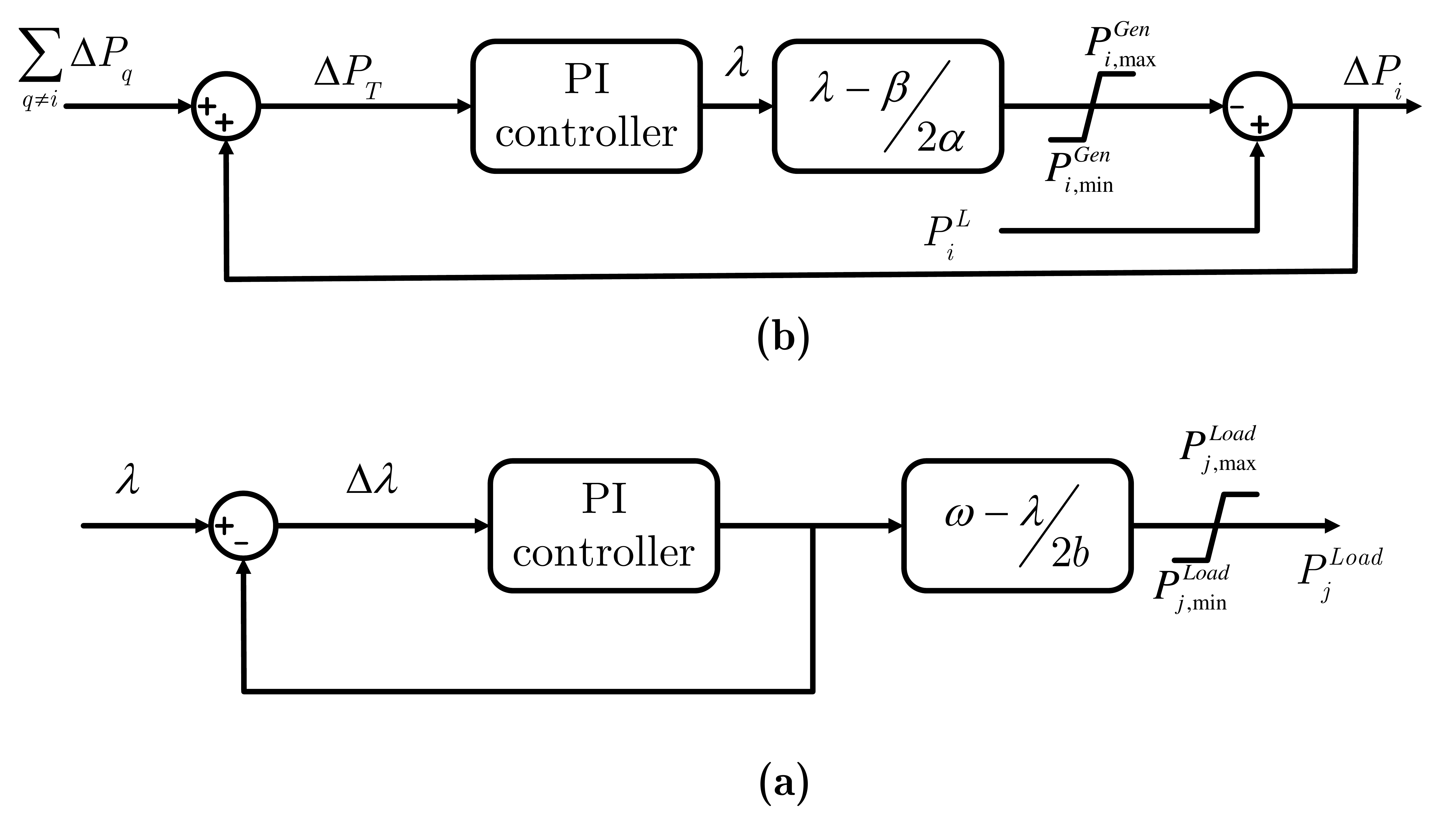}

\caption{Local PI-controller for DGs and consumers; (a) PI-controller for DGs,
(b) PI-controller for consumers\label{fig:PIcontroller}}
\end{figure}

The proposed distributed algorithm is particularly designed for easy
implementation and configuration of local agents by using a simple
proportional-integral (PI) or integral (I) controller. Equations (\ref{eq:completedeltaP})-(\ref{eq:powerimbalance})
can be easily modeled as a PI or I controller to update the estimated
power mismatch iteratively and exchange information (\textit{i.e.},
$\Delta P$) with other agents. Figure \ref{fig:PIcontroller}(a)
shows a simple PI controller for a DG. A PI controller is used for
each consumer to adjust its own demand based on $\Delta\lambda$ between
the current and previous iteration. Figure \ref{fig:PIcontroller}(b)
shows a simple controller for a consumer. As $\Delta\lambda$ becomes
zero, the consumer\textquoteright s surplus is maximized.

\begin{figure}[tbph] \centering
\includegraphics[width=0.7\columnwidth]{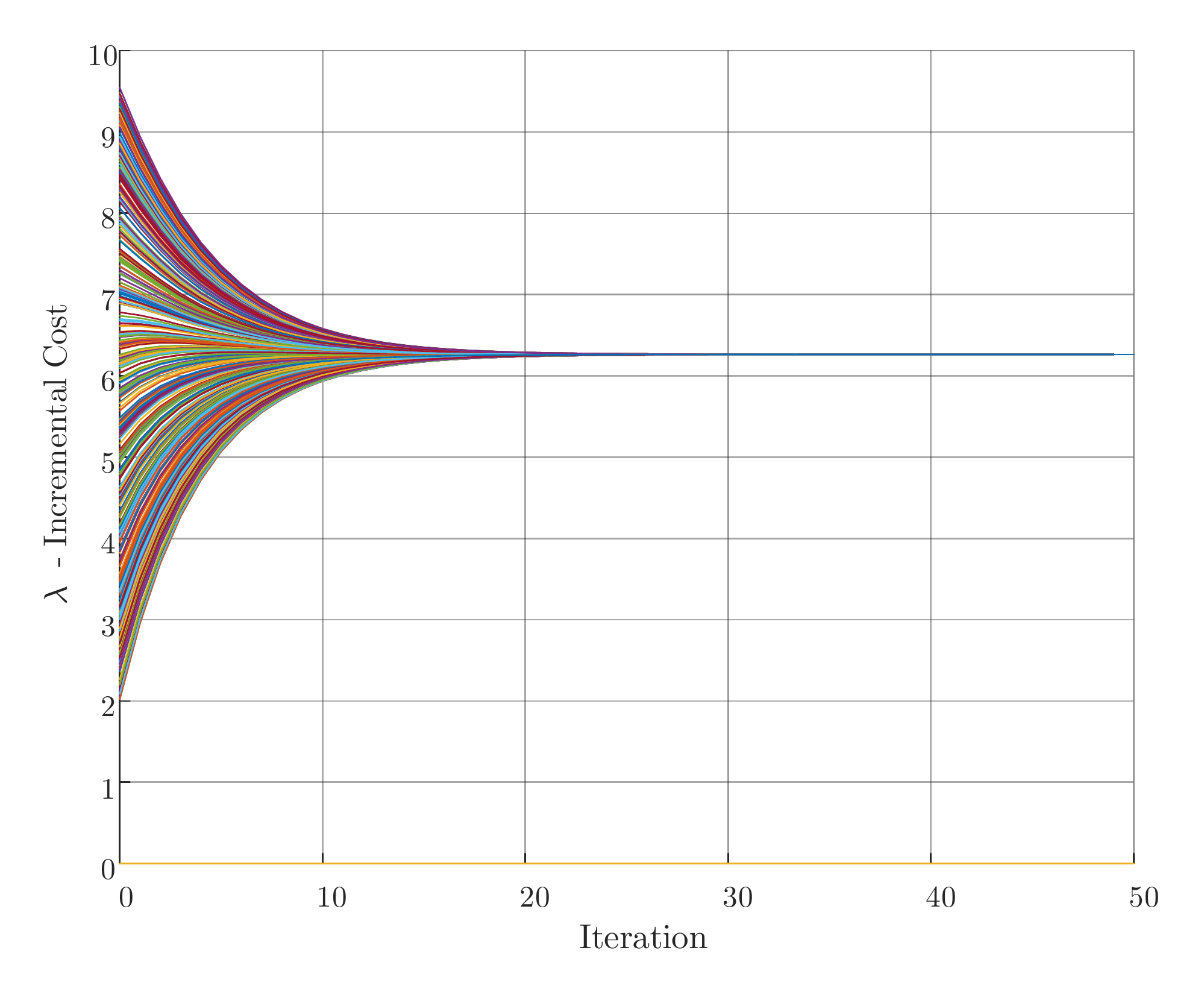}

\caption{Incremental cost ($\unitfrac{\$}{kWh}$) of DGs in case study II\label{fig:1400v2}}
\end{figure}

\begin{figure}[tbph] \centering
\includegraphics[width=0.7\columnwidth]{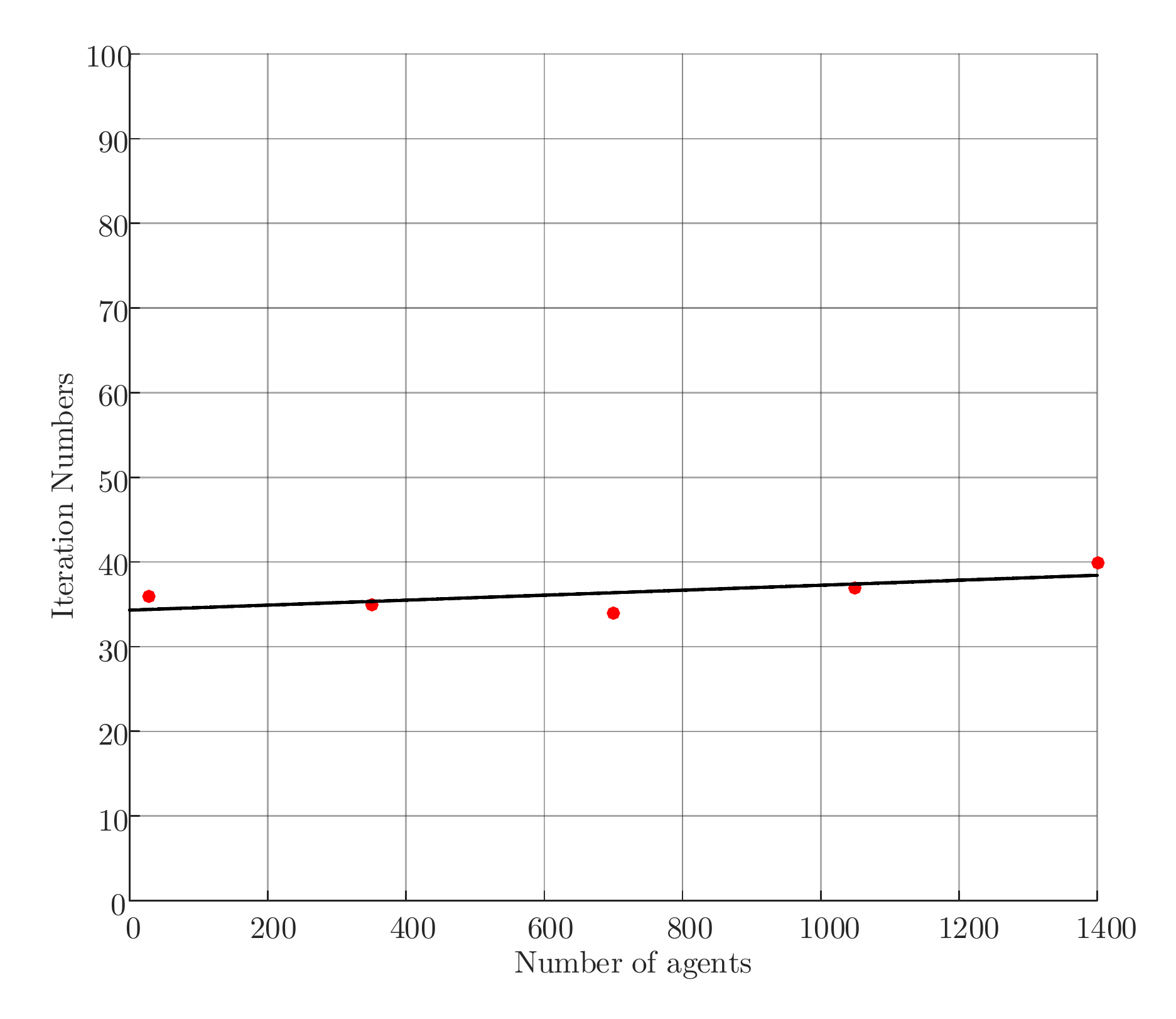}

\caption{The trend in number of iteration for algorithm conv\label{fig:Thetrend}}
\end{figure}

In case study III, the communication platform is
implemented in VOLTTRON\texttrademark{} which is an innovative distributed
control and sensing software platform developed by the Pacific Northwest
National Laboratory \cite{Luo2017}. The open-source VOLTTRON\texttrademark{}
platform makes it possible to deploy distributed control agents at
a very low cost. The platform provides various services such as resource
management, agent code verification and directory services allowing
to manage different assets within the power system. In the large scale,
VOLTTRON\texttrademark{} can manage assets within smart grids, facilitate
demand response, support energy trading and record grid data.

\begin{figure}[tbh] \centering
\includegraphics[width=0.5\columnwidth]{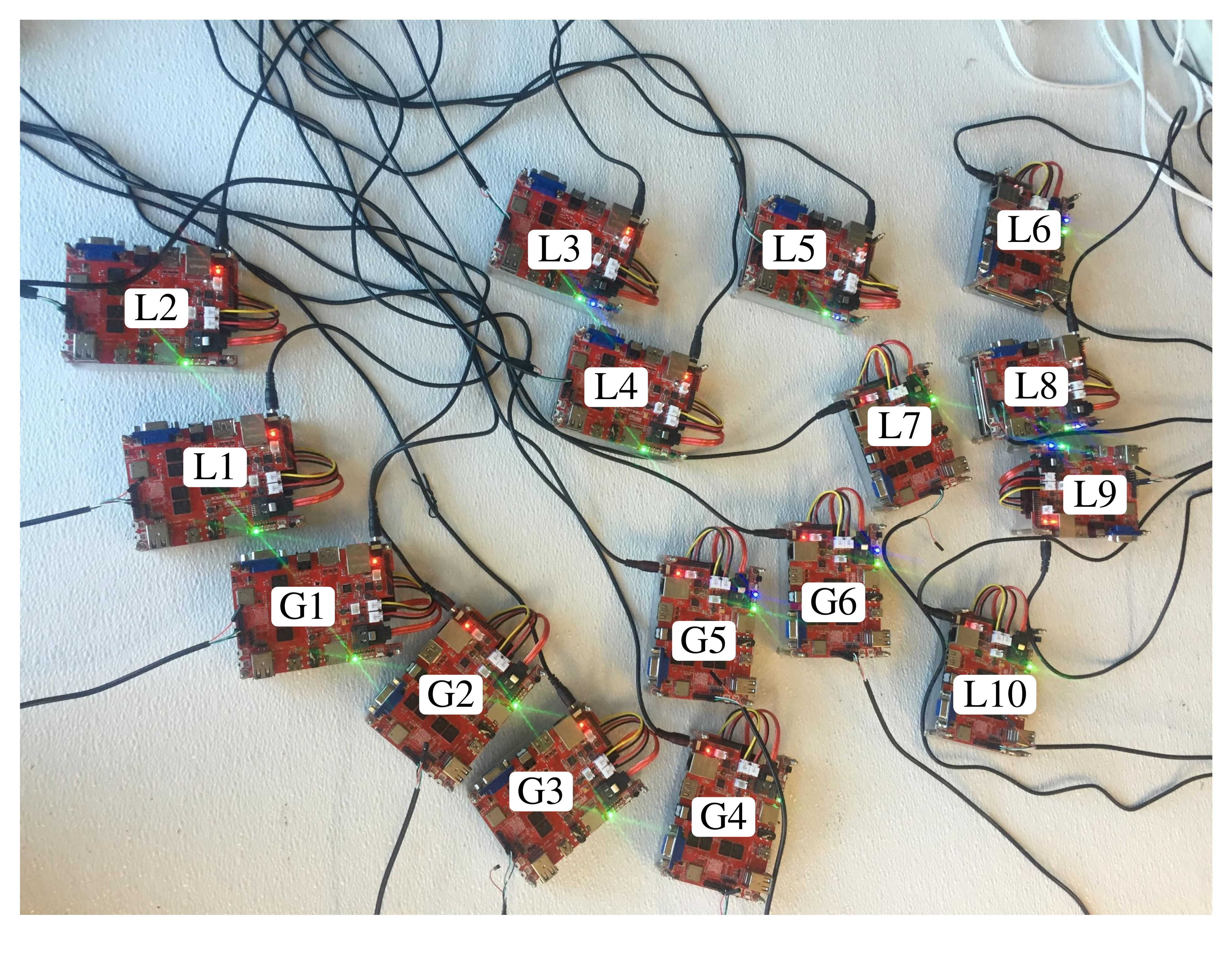}

\caption{The Proof-of-Concept 16-node Testbed using a Cluster of Single-board
PCs (``L'' inidicates loads/consumers; ``G'' inididates DGs)\label{fig:OverallTestbedArchitecture}}
\end{figure}

The VOLTTRON\texttrademark{} platform is implemented
on an ordinary Linux desktop with an FX-4100 CPU @ 3.6 GHz, 8-GB RAM
memory. The control platform is substantiated into
a cluster of low-cost credit-card-size single board PCs (Cubieboard
A20). The Cubieboard A20 processor is based on
a dual-core ARM Cortec-A7 CPU architecture. We use the Python programming
language to implement the proposed consensus-based distributed control
algorithms for each agent. In this proof-of-concept implementation,
every Cubieboard is emulated as a distributed controller for local
agents (DGs and consumers), while the PC with {} the
VOLTTRON\texttrademark{} platform is regarded as an information exchange
bus. The decision making process is conducted in a fully distributed
fashion. Figure \ref{fig:OverallTestbedArchitecture} shows the overall
system set up. For the demonstration purpose, the proposed distributed
algorithm is applied to a relatively small-scale distribution system
including 6 DGs and 10 consumers. As shown in figure \ref{fig:OverallTestbedArchitecture},
DGs are labeled as G1, G2, ..., G6 while consumers are labeled as
L1, L2, ..., L10. The coefficients of DGs' cost functions and consumers'
utility functions are randomly selected.

\begin{figure}[H] \centering
\center\includegraphics[width=0.5\columnwidth]{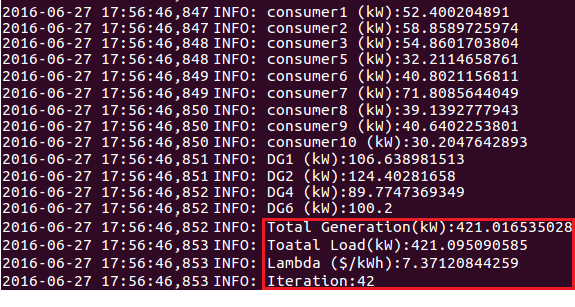}

\caption{The experimental test results\label{fig:Experimentalresults}}
\end{figure}

Figure \ref{fig:Experimentalresults} shows the detailed experimental
test results. The output-power of DG3 and DG5 and the demand of L4
are zero. The experimental test results are validated using the benchmark
results achieved by a centralized approach. As figure \ref{fig:Experimentalresults}
shows, the total generation is about $\unit[421]{kW}$ that satisfies
the total load demand. The incremental cost converges to $\unitfrac[7.371]{\$}{kwh}${}
at the 42-th iteration, which is considered as a fast convergence
rate.

\FloatBarrier

\section{Conclusions\label{sec:ConclusionsandFuture}}

In this paper, we proposed a novel consensus-based distributed algorithm
to solve an optimal dispatch problem of distributed generators. First,
we formulated the social welfare problem considering the cost functions
of DGs and the utility functions of consumers. Second, we developed
the distributed algorithm to find the global optima by allowing the
iterative coordination of agents (consumers and DGs) with each other.
Agents only share their estimated power mismatch, which does not contain
any private information, ultimately contributing to a fair electricity
market. Third, we performed software simulation and experimental test
to demonstrate the accuracy, privacy, effectiveness, fast-convergence,
scalability, and easy-implementation of the proposed distributed algorithm
under various conditions.

\bibliographystyle{IEEEtran}
\bibliography{library}

\begin{thebibliography}{10}
\providecommand{\url}[1]{#1}
\csname url@samestyle\endcsname
\providecommand{\newblock}{\relax}
\providecommand{\bibinfo}[2]{#2}
\providecommand{\BIBentrySTDinterwordspacing}{\spaceskip=0pt\relax}
\providecommand{\BIBentryALTinterwordstretchfactor}{4}
\providecommand{\BIBentryALTinterwordspacing}{\spaceskip=\fontdimen2\font plus
\BIBentryALTinterwordstretchfactor\fontdimen3\font minus
  \fontdimen4\font\relax}
\providecommand{\BIBforeignlanguage}[2]{{%
\expandafter\ifx\csname l@#1\endcsname\relax
\typeout{** WARNING: IEEEtran.bst: No hyphenation pattern has been}%
\typeout{** loaded for the language `#1'. Using the pattern for}%
\typeout{** the default language instead.}%
\else
\language=\csname l@#1\endcsname
\fi
#2}}
\providecommand{\BIBdecl}{\relax}
\BIBdecl

\bibitem{Fang2012}
\BIBentryALTinterwordspacing
X.~Fang, S.~Misra, G.~Xue, and D.~Yang, ``{Smart grid - The new and improved
  power grid: A survey},'' \emph{IEEE Communications Surveys and Tutorials},
  vol.~14, no.~4, pp. 944--980, jan 2012. [Online]. Available:
  \url{http://ieeexplore.ieee.org/lpdocs/epic03/wrapper.htm?arnumber=6099519}
\BIBentrySTDinterwordspacing

\bibitem{Potter2009}
C.~W. Potter, A.~Archambault, and K.~Westrick, ``{Building a smarter smart grid
  through better renewable energy information},'' \emph{2009 IEEE/PES Power
  Systems Conference and Exposition, PSCE 2009}, pp. 1--5, 2009.

\bibitem{NationalAcademiesofSciencesandMedicine2016}
\BIBentryALTinterwordspacing
\emph{\BIBforeignlanguage{English}{{Analytic Research Foundations for the
  Next-Generation Electric Grid}}}.\hskip 1em plus 0.5em minus 0.4em\relax
  Washington, D.C.: National Academies Press, 2016. [Online]. Available:
  \url{http://www.nap.edu/catalog/21919}
\BIBentrySTDinterwordspacing

\bibitem{Asadinejad}
A.~Asadinejad, K.~Tomsovic, and C.~Chen, ``{Sensitivity of Incentive Based
  Demand Response Program To Residential Customer Elasticity},'' in \emph{North
  American Power Symposium(NAPS)}.\hskip 1em plus 0.5em minus 0.4em\relax
  Denver, CO: IEEE, sep 2016.

\bibitem{Mudumbai2012c}
\BIBentryALTinterwordspacing
R.~Mudumbai, S.~Dasgupta, and B.~B. Cho, ``{Distributed control for optimal
  economic dispatch of a network of heterogeneous power generators},''
  \emph{IEEE Transactions on Power Systems}, vol.~27, no.~4, pp. 1750--1760,
  nov 2012. [Online]. Available:
  \url{http://ieeexplore.ieee.org/lpdocs/epic03/wrapper.htm?arnumber=6178300}
\BIBentrySTDinterwordspacing

\bibitem{Saber2003}
\BIBentryALTinterwordspacing
R.~Saber and R.~Murray, ``{Consensus protocols for networks of dynamic
  agents},'' in \emph{Proceedings of the 2003 American Control Conference,
  2003.}, vol.~2.\hskip 1em plus 0.5em minus 0.4em\relax IEEE, 2003, pp.
  951--956. [Online]. Available:
  \url{http://ieeexplore.ieee.org/lpdocs/epic03/wrapper.htm?arnumber=1239709}
\BIBentrySTDinterwordspacing

\bibitem{Ren}
\BIBentryALTinterwordspacing
W.~Ren and R.~W.Beard, \emph{{Distributed consensus in multi-vechicle
  cooperative control}}, ser. Communications and Control Engineering.\hskip 1em
  plus 0.5em minus 0.4em\relax London: Springer London, 2008, vol.~36.
  [Online]. Available: \url{http://www.ncbi.nlm.nih.gov/pubmed/24815723}
\BIBentrySTDinterwordspacing

\bibitem{Olfati-Saber2007}
\BIBentryALTinterwordspacing
R.~Olfati-Saber, J.~A. Fax, and R.~M. Murray, ``{Consensus and cooperation in
  networked multi-agent systems},'' \emph{Proceedings of the IEEE}, vol.~95,
  no.~1, pp. 215--233, jan 2007. [Online]. Available:
  \url{http://ieeexplore.ieee.org/lpdocs/epic03/wrapper.htm?arnumber=4118472}
\BIBentrySTDinterwordspacing

\bibitem{Binetti2014}
\BIBentryALTinterwordspacing
G.~Binetti, A.~Davoudi, F.~L. Lewis, D.~Naso, and B.~Turchiano, ``{Distributed
  consensus-based economic dispatch with transmission losses},'' \emph{IEEE
  Transactions on Power Systems}, vol.~29, no.~4, pp. 1711--1720, 2014.
  [Online]. Available:
  \url{http://ieeexplore.ieee.org/lpdocs/epic03/wrapper.htm?arnumber=6717171}
\BIBentrySTDinterwordspacing

\bibitem{Xu2017}
\BIBentryALTinterwordspacing
S.~Xu, H.~Pourbabak, and W.~Su, ``{Distributed cooperative control for economic
  operation of multiple plug-in electric vehicle parking decks},''
  \emph{International Transactions on Electrical Energy Systems}, p. e2348,
  2017. [Online]. Available: \url{http://doi.wiley.com/10.1002/etep.2348}
\BIBentrySTDinterwordspacing

\bibitem{Xu2016Su}
\BIBentryALTinterwordspacing
Y.~Xu, J.~Hu, W.~Gu, W.~Su, and W.~Liu, ``{Real-Time Distributed Control of
  Battery Energy Storage Systems for Security Constrained DC-OPF},'' \emph{IEEE
  Transactions on Smart Grid}, pp. 1--1, 2016. [Online]. Available:
  \url{http://ieeexplore.ieee.org/document/7523948/}
\BIBentrySTDinterwordspacing

\bibitem{Zhang2012d}
\BIBentryALTinterwordspacing
Z.~Zhang and M.~Y. Chow, ``{Convergence analysis of the incremental cost
  consensus algorithm under different communication network topologies in a
  smart grid},'' \emph{IEEE Transactions on Power Systems}, vol.~27, no.~4, pp.
  1761--1768, nov 2012. [Online]. Available:
  \url{http://ieeexplore.ieee.org/lpdocs/epic03/wrapper.htm?arnumber=6183499}
\BIBentrySTDinterwordspacing

\bibitem{Pourbabakbook}
\BIBentryALTinterwordspacing
H.~Pourbabak, T.~Chen, B.~Zhang, and W.~Su, ``{Control and energy management
  system in microgrids},'' in \emph{Clean Energy Microgrids}.\hskip 1em plus
  0.5em minus 0.4em\relax Institution of Engineering and Technology, 2017,
  ch.~3, pp. 109--133. [Online]. Available:
  \url{http://digital-library.theiet.org/content/books/10.1049/pbpo090e_ch3}
\BIBentrySTDinterwordspacing

\bibitem{Zeng2014}
\BIBentryALTinterwordspacing
W.~Zeng and M.~Y. Chow, ``{Resilient distributed control in the presence of
  misbehaving agents in networked control systems},'' \emph{IEEE Transactions
  on Cybernetics}, vol.~44, no.~11, pp. 2038--2049, nov 2014. [Online].
  Available: \url{http://www.ncbi.nlm.nih.gov/pubmed/25330469}
\BIBentrySTDinterwordspacing

\bibitem{Kirschen2004}
\BIBentryALTinterwordspacing
D.~Kirschen and G.~Strbac, \emph{{Fundamentals of Power System
  Economics}}.\hskip 1em plus 0.5em minus 0.4em\relax Chichester, UK: John
  Wiley {\&} Sons, Ltd, mar 2004, vol.~4, no.~4. [Online]. Available:
  \url{http://doi.wiley.com/10.1002/0470020598}
\BIBentrySTDinterwordspacing

\bibitem{Samadi2012}
\BIBentryALTinterwordspacing
P.~Samadi, H.~Mohsenian-Rad, R.~Schober, and V.~W.~S. Wong, ``{Advanced demand
  side management for the future smart grid using mechanism design},''
  \emph{IEEE Transactions on Smart Grid}, vol.~3, no.~3, pp. 1170--1180, sep
  2012. [Online]. Available:
  \url{http://ieeexplore.ieee.org/lpdocs/epic03/wrapper.htm?arnumber=6266724}
\BIBentrySTDinterwordspacing

\bibitem{Xu2015b}
\BIBentryALTinterwordspacing
Y.~Xu, ``{Optimal Distributed Charging Rate Control of Plug-In Electric
  Vehicles for Demand Management},'' \emph{IEEE Transactions on Power Systems},
  vol.~30, no.~3, pp. 1536--1545, may 2015. [Online]. Available:
  \url{http://ieeexplore.ieee.org/document/6889041/
  http://ieeexplore.ieee.org/lpdocs/epic03/wrapper.htm?arnumber=6889041}
\BIBentrySTDinterwordspacing

\bibitem{Soediono1989}
\BIBentryALTinterwordspacing
W.~Ren and Y.~Cao, \emph{{Distributed Coordination of Multi-agent Networks}},
  ser. Communications and Control Engineering, {Intergovernmental Panel on
  Climate Change}, Ed.\hskip 1em plus 0.5em minus 0.4em\relax London: Springer
  London, nov 2011, vol.~58, no.~12. [Online]. Available:
  \url{http://link.springer.com/10.1007/978-0-85729-169-1}
\BIBentrySTDinterwordspacing

\bibitem{Kazemi2016}
\BIBentryALTinterwordspacing
A.~Kazemi and H.~Pourbabak, ``{Islanding detection method based on a new
  approach to voltage phase angle of constant power inverters},'' \emph{IET
  Generation, Transmission {\&} Distribution}, vol.~10, no.~5, pp. 1190--1198,
  apr 2016. [Online]. Available:
  \url{http://digital-library.theiet.org/content/journals/10.1049/iet-gtd.2015.0776}
\BIBentrySTDinterwordspacing

\bibitem{Guo2014}
\BIBentryALTinterwordspacing
F.~Guo, C.~Wen, J.~Mao, and Y.~D. Song, ``{Distributed Economic Dispatch for
  Smart Grids with Random Wind Power},'' \emph{IEEE Transactions on Smart
  Grid}, vol.~7, no.~3, pp. 1572--1583, may 2016. [Online]. Available:
  \url{http://ieeexplore.ieee.org/lpdocs/epic03/wrapper.htm?arnumber=7120161}
\BIBentrySTDinterwordspacing

\bibitem{Yi2015}
\BIBentryALTinterwordspacing
P.~Yi, Y.~Hong, and F.~Liu, ``{Distributed gradient algorithm for constrained
  optimization with application to load sharing in power systems},''
  \emph{Systems {\&} Control Letters}, vol.~83, pp. 45--52, sep 2015. [Online].
  Available:
  \url{http://www.sciencedirect.com/science/article/pii/S0167691115001346}
\BIBentrySTDinterwordspacing

\bibitem{Hug2015a}
\BIBentryALTinterwordspacing
G.~Hug, S.~Kar, and C.~Wu, ``{Consensus + Innovations Approach for Distributed
  Multiagent Coordination in a Microgrid},'' \emph{IEEE Transactions on Smart
  Grid}, vol.~6, no.~4, pp. 1893--1903, jul 2015. [Online]. Available:
  \url{http://ieeexplore.ieee.org/lpdocs/epic03/wrapper.htm?arnumber=7098427}
\BIBentrySTDinterwordspacing

\bibitem{Rahbari-Asr2014a}
\BIBentryALTinterwordspacing
N.~Rahbari-Asr, U.~Ojha, Z.~Zhang, and M.~Y. Chow, ``{Incremental welfare
  consensus algorithm for cooperative distributed generation/demand response in
  smart grid},'' \emph{IEEE Transactions on Smart Grid}, vol.~5, no.~6, pp.
  2836--2845, nov 2014. [Online]. Available:
  \url{http://ieeexplore.ieee.org/lpdocs/epic03/wrapper.htm?arnumber=6884808}
\BIBentrySTDinterwordspacing

\bibitem{Rahbari-Asr2014}
\BIBentryALTinterwordspacing
N.~Rahbari-Asr and M.~Y. Chow, ``{Cooperative distributed demand management for
  community charging of PHEV/PEVs based on KKT conditions and consensus
  networks},'' \emph{IEEE Transactions on Industrial Informatics}, vol.~10,
  no.~3, pp. 1907--1916, 2014. [Online]. Available:
  \url{http://ieeexplore.ieee.org/lpdocs/epic03/wrapper.htm?arnumber=6730946}
\BIBentrySTDinterwordspacing

\bibitem{DallAnese2013}
E.~Dall'Anese, H.~Zhu, and G.~B. Giannakis, ``{Distributed optimal power flow
  for smart microgrids},'' \emph{IEEE Transactions on Smart Grid}, vol.~4,
  no.~3, pp. 1464--1475, 2013.

\bibitem{Multibuyer2015}
\BIBentryALTinterwordspacing
R.~Deng, Z.~Yang, F.~Hou, M.~Chow, and J.~Chen, ``{Distributed Real Time Demand
  Response in Multiseller Multibuyer Smart Distribution Grid},'' \emph{IEEE
  Transactions on Power Systems}, vol.~30, no.~5, pp. 2364--2374, sep 2015.
  [Online]. Available: \url{http://ieeexplore.ieee.org/document/6922169/}
\BIBentrySTDinterwordspacing

\bibitem{Elsayed2015}
\BIBentryALTinterwordspacing
W.~T. Elsayed and E.~F. El-Saadany, ``{A Fully Decentralized Approach for
  Solving the Economic Dispatch Problem},'' \emph{IEEE Transactions on Power
  Systems}, vol.~30, no.~4, pp. 2179--2189, jul 2014. [Online]. Available:
  \url{http://ieeexplore.ieee.org/lpdocs/epic03/wrapper.htm?arnumber=6917059}
\BIBentrySTDinterwordspacing

\bibitem{Zhang2015b}
\BIBentryALTinterwordspacing
W.~Zhang, W.~Liu, X.~Wang, L.~Liu, and F.~Ferrese, ``{Online optimal generation
  control based on constrained distributed gradient algorithm},'' \emph{IEEE
  Transactions on Power Systems}, vol.~30, no.~1, pp. 35--45, jan 2015.
  [Online]. Available:
  \url{http://ieeexplore.ieee.org/lpdocs/epic03/wrapper.htm?arnumber=6810888}
\BIBentrySTDinterwordspacing

\bibitem{Yang2013a}
\BIBentryALTinterwordspacing
S.~Yang, S.~Tan, and J.-X. Xu, ``{Consensus Based Approach for Economic
  Dispatch Problem in a Smart Grid},'' \emph{IEEE Transactions on Power
  Systems}, vol.~28, no.~4, pp. 4416--4426, 2013. [Online]. Available:
  \url{http://ieeexplore.ieee.org/lpdocs/epic03/wrapper.htm?arnumber=6560423}
\BIBentrySTDinterwordspacing

\bibitem{Cai2012}
\BIBentryALTinterwordspacing
N.~Cai, N.~T.~T. Nga, and J.~Mitra, ``{Economic dispatch in microgrids using
  multi-agent system},'' in \emph{2012 North American Power Symposium, NAPS
  2012}.\hskip 1em plus 0.5em minus 0.4em\relax IEEE, sep 2012, pp. 1--5.
  [Online]. Available:
  \url{http://ieeexplore.ieee.org/lpdocs/epic03/wrapper.htm?arnumber=6336435}
\BIBentrySTDinterwordspacing

\bibitem{Binetti2014b}
G.~Binetti, A.~Davoudi, D.~Naso, B.~Turchiano, and F.~L. Lewis, ``{A
  distributed auction-based algorithm for the nonconvex economic dispatch
  problem},'' \emph{IEEE Transactions on Industrial Informatics}, vol.~10,
  no.~2, pp. 1124--1132, 2014.

\bibitem{Chang2014}
\BIBentryALTinterwordspacing
T.-H. Chang, A.~Nedic, and A.~Scaglione, ``{Distributed Constrained
  Optimization by Consensus-Based Primal-Dual Perturbation Method},''
  \emph{IEEE Transactions on Automatic Control}, vol.~59, no.~6, pp.
  1524--1538, jun 2014. [Online]. Available:
  \url{http://ieeexplore.ieee.org/document/6748910/}
\BIBentrySTDinterwordspacing

\bibitem{Pourbabakconf}
\BIBentryALTinterwordspacing
H.~Pourbabak, {Tao Chen}, and W.~Su, ``{Consensus-based distributed control for
  economic operation of distribution grid with multiple consumers and
  prosumers},'' in \emph{2016 IEEE Power and Energy Society General Meeting
  (PESGM)}, IEEE.\hskip 1em plus 0.5em minus 0.4em\relax Boston, MA, July
  17-21, 2016: IEEE, jul 2016, pp. 1--5. [Online]. Available:
  \url{http://ieeexplore.ieee.org/document/7741083/}
\BIBentrySTDinterwordspacing

\bibitem{Fahrioglu2001}
\BIBentryALTinterwordspacing
M.~Fahrioglu and F.~Alvarado, ``{Using utility information to calibrate
  customer demand management behavior models},'' \emph{IEEE Transactions on
  Power Systems}, vol.~16, no.~2, pp. 317--322, may 2001. [Online]. Available:
  \url{http://ieeexplore.ieee.org/lpdocs/epic03/wrapper.htm?arnumber=918305
  http://ieeexplore.ieee.org/document/918305/}
\BIBentrySTDinterwordspacing

\bibitem{Samadi2010}
\BIBentryALTinterwordspacing
P.~Samadi, A.-H. Mohsenian-Rad, R.~Schober, V.~W.~S. Wong, and J.~Jatskevich,
  ``{Optimal Real-Time Pricing Algorithm Based on Utility Maximization for
  Smart Grid},'' \emph{2010 First IEEE International Conference on Smart Grid
  Communications}, pp. 415--420, 2010. [Online]. Available:
  \url{http://ieeexplore.ieee.org/lpdocs/epic03/wrapper.htm?arnumber=5622077}
\BIBentrySTDinterwordspacing

\bibitem{Mohajeryami2016}
\BIBentryALTinterwordspacing
S.~Mohajeryami, I.~N. Moghaddam, M.~Doostan, B.~Vatani, and P.~Schwarz, ``{A
  novel economic model for price-based demand response},'' \emph{Electric Power
  Systems Research}, vol. 135, pp. 1--9, jun 2016. [Online]. Available:
  \url{http://linkinghub.elsevier.com/retrieve/pii/S0378779616300736}
\BIBentrySTDinterwordspacing

\bibitem{Nechyba2016}
\BIBentryALTinterwordspacing
T.~J. Nechyba, \emph{{Microeconomics: An Intuitive Approach With
  Calculus}}.\hskip 1em plus 0.5em minus 0.4em\relax Cengage Learning, 2010.
  [Online]. Available:
  \url{http://books.google.com/books?id=T4WuH3PRtYkC\&pgis=1}
\BIBentrySTDinterwordspacing

\bibitem{Melorose2015a}
\BIBentryALTinterwordspacing
F.~L. Lewis, H.~Zhang, K.~Hengster-Movric, and A.~Das, \emph{{Cooperative
  Control of Multi-Agent Systems}}, ser. Communications and Control
  Engineering.\hskip 1em plus 0.5em minus 0.4em\relax London: Springer London,
  2014, vol. 1542. [Online]. Available:
  \url{http://link.springer.com/10.1007/978-1-4471-5574-4}
\BIBentrySTDinterwordspacing

\bibitem{Boyd2004}
S.~P. Boyd and L.~Vandenberghe, \emph{{Convex optimization}}.\hskip 1em plus
  0.5em minus 0.4em\relax Cambridge UK ;New York: Cambridge University Press,
  2004.

\bibitem{Ruszczynski2011}
A.~Ruszczynski, \emph{{Nonlinear Optimization}}.\hskip 1em plus 0.5em minus
  0.4em\relax Princeton University Press, 2006.

\bibitem{YALMIPTOOLBOX}
\BIBentryALTinterwordspacing
J.~Lofberg, ``{YALMIP : a toolbox for modeling and optimization in MATLAB},''
  in \emph{2004 IEEE International Conference on Computer Aided Control Systems
  Design}.\hskip 1em plus 0.5em minus 0.4em\relax IEEE, 2004, pp. 284--289.
  [Online]. Available:
  \url{http://ieeexplore.ieee.org/lpdocs/epic03/wrapper.htm?arnumber=1393890}
\BIBentrySTDinterwordspacing

\bibitem{Luo2017}
\BIBentryALTinterwordspacing
{Jingwei Luo}, H.~Pourbabak, and W.~Su, ``{The application of distributed
  control algorithms using VOLTTRON-based software platform},'' in \emph{2017
  8th International Renewable Energy Congress (IREC)}.\hskip 1em plus 0.5em
  minus 0.4em\relax Dead Sea, Jordan: IEEE, mar 2017, pp. 1--6. [Online].
  Available: \url{http://ieeexplore.ieee.org/document/7926056/}
\BIBentrySTDinterwordspacing

\end{thebibliography}

\end{document}